\documentclass[10pt]{article}

\usepackage{a4wide}
\usepackage{amssymb}
\usepackage{amsfonts}
\usepackage{amsmath}
\input xy
\xyoption{arrow} \xyoption{matrix}

\date{}

\newtheorem{proposition}{Proposition}[section]
\newtheorem{theorem}[proposition]{Theorem}
\newtheorem{lemma}[proposition]{Lemma}

\newtheorem{corollary}[proposition]{Corollary}

\def\der{\partial }

\def\nFM0{{\nu }_{F,M_0}}
\def\nFN0{{\nu }_{F,N_0}}
\def\nGN0{{\nu }_{G,N_0}}

\def\N0{ {\bf N}_0 }

\def\t{\otimes}
\def\g{\gamma}
\def\v{\varphi}
\def\ra{\rightarrow}

\def\lra{\leftrightarrow}
\def\Xpm{X^{\pm }}

\def\Z{\mathbb{Z}}

\def\l1{{\lambda}_1}

\def\a{\alpha}
\def\a0{ {\alpha }_0}
\def\a1{ {\alpha }_1}

\def\l{\lambda}
\def\o{\omega}

\def\nFGM0{{\nu }_{F,G,M_0}}

\def\nFN0{{\nu}_{F,N_0}}

\def\sm{{\sigma}^m}

\def\sm1{{\sigma}^{-1}}

\def\smtp1{{\sigma}^{-t+1}}

\def\o{\omega }
\def\S1{S^{-1}}

\def\Xpm1{X^{\pm 1}_1}

\def\sPM1{{\sigma }^{\pm 1}}
\def\sMP1{{\sigma }^{\mp 1 }}

\def\d{\delta}

\def\di{{\rm d.ind}}

\def\L{\Lambda}

\def\G{\Gamma}

\def\Ytm1{Y^{t-1}}
\def\Yim1{Y^{i-1}}

\def\CS{{\cal S}}

\def\CH{{\cal H}}

\def\CZ{{\cal Z}}
\def\supp{{\rm supp }}

\def\Aut{{\rm Aut}}

\def\char{{\rm char }}

\def\ker{ {\rm ker } }

\def\CJ{ {\cal J}}

\def\D{ \Delta }

\def\SL2Z{ {\rm SL}_2({\bf Z}) }

\def\CZ{ {\cal Z}}

\def\th{ \theta }

\def\Gp1{ G^{1 , 1 } }
\def\P11{ P^{-1 , 1 } }
\def\Pp1{ P^{1 , 1 } }

\def\th{\theta}
\def\CE{{\cal E}}

\def\nCLsr{{}^\nu\kern-2pt {\cal L}^{\sigma , \rho  }}
\def\nP{{}^\nu \kern-2pt P}
\def\nL{{}^\nu\kern-2pt L}
\def\nLL{{}^\nu\kern-2pt \Lambda}
\def\nPsr{{}^\nu\kern-2pt P^{\sigma , \rho  }}
\def\nLsr{{}^\nu\kern-2pt L^{\sigma , \rho  }}
\def\nuCL{{}^\nu\kern-2pt  {\cal L}}
\def\nCLsr{{}^\nu\kern-2pt {\cal L}^{\sigma , \rho  }}
\def\nCL1m{{}^\nu\kern-2pt {\cal L}^{-1 , 1  }}
\def\x1nu{x^\frac{1}{\nu}}
\def\xm1nu{x^{-\frac{1}{\nu}}}

\def\ra{\rightarrow }

\def\CE{ {\cal E} }
\def\CH{ {\cal H}}

\def\nAM0{{\nu }_{{\cal A},M_0}}
\def\nAN0{{\nu }_{{\cal A},N_0}}

\def\End{ {\rm End }}

\def\CJ{ {\cal J }}

\def\det{ {\rm det }}

\def\ga{\mathfrak{a}}
\def\gb{\mathfrak{b}}

\def\gp{\mathfrak{p}}

\def\GL{{\rm GL}}
\def\SL{{\rm SL}}

\def\di!{\frac{\der^i}{i!}}
\def\dik!{\frac{\der^k_i}{k!}}

\def\N{\mathbb{N}}

\def\0{\overline{0}}
\def\1{\overline{1}}

\def\Ln1{\L_{n,\overline{1}}}

\def\a1{a_{\overline{1}}}

\def\S{\Sigma}

\def\vn1{\overrightarrow{n-1}}

\def\im{{\rm im}}

\def\mA{\mathbb{A}}

\def\bu{\overline{u}}

\def\Inn{{\rm Inn}}

\def\mS{\mathbb{S}}

\def\clKdim{{\rm cl.Kdim}}

\def\mT{\mathbb{T}}

\def\mE{\mathbb{E}}
\def\mX{\mathbb{X}}

\def\bgp{\overline{\gp}}

\def\K1{{\rm K}_1}
\def\bG{\overline{\G}}
\def\bpsi{\overline{\psi}}
\def\mK{\mathbb{K}}
\def\mP{\mathbb{P}}
\def\mY{\mathbb{Y}}

\def\tTh{\widetilde{\Theta}}
\def\tmE{\widetilde{\mathbb{E}}}
\def\bdet{\overline{\det}}

\begin{document}

\author{V. V. \  Bavula 
}

\title{The group ${\rm K}_1(\mS_n)$  of the  algebra
 of one-sided inverses of a polynomial algebra}

\maketitle

\begin{abstract}
The algebra $\mS_n$ of one-sided inverses of a polynomial algebra
$P_n$ in $n$ variables is obtained from $P_n$ by adding commuting,
 {\em left} (but not two-sided) inverses of the canonical
 generators of the algebra $P_n$. The algebra $\mS_n$ is a
 noncommutative, non-Noetherian algebra of classical Krull
 dimension $2n$ and of global dimension $n$ and is not a domain. If the ground field $K$ has
 characteristic zero then the algebra $\mS_n$ is canonically
 isomorphic to the algebra $K\langle
\frac{\der}{\der x_1}, \ldots ,\frac{\der}{\der x_n},  \int_1,
\ldots , \int_n\rangle $ of  scalar  integro-differential
operators. 
  It is proved that $ {\rm K}_1(\mS_n)\simeq K^*$. The main idea is  to
 show that the group $\GL_\infty (\mS_n)$ is generated by $K^*$,
  the group of elementary matrices $E_\infty (\mS_n)$ and  $(n-2)2^{n-1}+1$ explicit (tricky)
 matrices and then to prove that all the matrices are
 elementary. For each nonzero idempotent prime ideal
 $\gp$ of height $m$ of the algebra $\mS_n$, it is proved that
 $$ {\rm K}_1(\mS_n, \gp )\simeq \begin{cases}
K^*, & \text{if }m=1,\\
\Z^{\frac{m(m-1)}{2}}\times K^{*m}& \text{if }m> 1.\\
\end{cases}$$




\end{abstract}


\section{Introduction}
 Throughout, ring
means an associative ring with $1$; module means a left module;
 $\N :=\{0, 1, \ldots \}$ is the set of natural numbers; $K$ is a
field  and  $K^*$ is its group of units; $P_n:= K[x_1, \ldots ,
x_n]$ is a polynomial algebra over $K$; $\der_1:=\frac{\der}{\der
x_1}, \ldots , \der_n:=\frac{\der}{\der x_n}$ are the usual partial
derivatives ($K$-linear derivations) of $P_n$; $\End_K(P_n)$ is
the algebra of all $K$-linear maps from $P_n$ to $P_n$ and
$\Aut_K(P_n)$ is its group of units (i.e. the group of all the
invertible linear maps from $P_n$ to $P_n$); the subalgebra $A_n:=
K \langle x_1, \ldots , x_n , \der_1, \ldots , \der_n\rangle$ of
$\End_K(P_n)$ is called the $n$'th {\em Weyl} algebra.

$\noindent $

{\it Definition}, \cite{shrekalg}. The 
{\em algebra} $\mathbb{S}_n=\mS_n(K)$ {\em of one-sided inverses}
of $P_n$ is the algebra generated over a field (or a ring) $K$ of
 by $2n$ elements $x_1, \ldots , x_n, y_n,
\ldots , y_n$ that satisfy the defining relations:
$$ y_1x_1=1, \ldots , y_nx_n=1, \;\; [x_i, y_j]=[x_i, x_j]= [y_i,y_j]=0
\;\; {\rm for\; all}\; i\neq j,$$ where $[a,b]:= ab-ba$ is  the
algebra commutator of elements $a$ and $b$.

$\noindent $

By the very definition, the algebra $\mS_n$ is obtained from the
polynomial algebra $P_n$ by adding commuting, left (but not
two-sided) inverses of its canonical generators. The algebra
$\mS_1=K\langle x, y\, | \, yx=1\rangle$ is a well-known primitive
algebra \cite{Jacobson-StrRing}, p. 35, Example 2.  Over the field
 $\mathbb{C}$ of complex numbers, the completion of the algebra
 $\mS_1$ is the {\em Toeplitz algebra} which is the
 $C^*$-algebra generated by a unilateral shift on the
 Hilbert space $l^2(\N )$ (note that $y=x^*$). The Toeplitz
 algebra is the universal $C^*$-algebra generated by a
 proper isometry. If $\char (K)=0$ then the algebra $\mS_n$ is
 isomorphic to the algebra $K\langle
\frac{\der}{\der x_1}, \ldots ,\frac{\der}{\der x_n},  \int_1,
\ldots , \int_n\rangle $ of  scalar  integro-differential
operators (via $x_i \mapsto \int_i$, $y_i\mapsto \frac{\der}{\der
x_i}$).

$\noindent $

In \cite{K1aut}, it is proved that ${\rm K}_1(\mS_1) \simeq K^*$.
The first aim of the paper is to prove that

\begin{itemize}
\item (Theorem \ref{17Apr10})  ${\rm K}_1(\mS_n) \simeq K^*$ for
all $n\geq 1$.
\end{itemize}

The algebra $\mS_n$ was studied in detail in \cite{shrekalg}: its
Gelfand-Kirillov dimension  is $2n$, its classical Krull dimension
$\clKdim (\mS_n)=2n$, and its   weak and  global dimensions  are $n$. The algebra $\mS_n$ is neither left nor right
Noetherian as was shown by Jacobson
\cite{Jacobson-onesidedinv-1950} when $n=1$ (see also Baer
\cite{Baer-inverses-1942}). Moreover, it contains infinite direct
sums of left and right ideals. It is an experimental fact that the
algebra $\mS_n \simeq \mS_1^{\t n}$ has properties that are
 a mixture of properties of the Weyl algebra $A_n\simeq A_1^{\t n}$
(in characteristic zero) and the polynomial algebra $P_{2n}\simeq
P_2^{\t n}$ which is not surprising when we look at their defining
relations:
\begin{eqnarray*}
 P_2=K\langle x,y\rangle :& yx-xy=0; \\
 A_1=K\langle x,y\rangle :&yx-xy=1;\\
 \mS_1=K\langle x,y\rangle :& yx=1.
\end{eqnarray*}
The group $G_n:= \Aut_{K-{\rm alg}}(\mS_n)$  of $K$-algebra
automorphisms of $\mS_n$ and the group $\mS_n^*$ of units of the algebra
$\mS_n$  were determined in the series of three papers \cite{shrekaut}, \cite{K1aut} and
\cite{Snaut},  and their explicit generators were found (both groups are
huge).  The group $G_1$ was found by Gerritzen
\cite{Gerritzen-yx=1-2000}).

\begin{theorem}\label{I1}
\begin{enumerate}
\item \cite{shrekaut} $G_n=S_n\ltimes \mT^n\ltimes \Inn (\mS_n)$
where $S_n$ is the symmetric group, $\mT^n\simeq K^{*n}$ is the
$n$-dimensional algebraic torus and  $\Inn (\mS_n)$ is the group
of inner automorphisms of $\mS_n$. \item \cite{K1aut},
\cite{jacaut} $\mS_n^* = K^*\times (1+\ga_n)^*$ where
$\ga_n$ is the ideal generated by  all the height one
prime ideals of $\mS_n$. \item  \cite{Snaut}  The
centre of the group $\mS_n^*$ is $K^*$ and the centre of the
group $(1+\ga_n)^*$ is $\{ 1 \}$. \item \cite{Snaut} The map
$(1+\ga_n)^*\ra \Inn (\mS_n)$, $u\mapsto \o_u$, is a group
isomorphism ($\o_u(a) = uau^{-1}$).
\end{enumerate}
\end{theorem}
{\bf The structure of the proof of Theorem \ref{17Apr10}}. The
idea of the proof that   ${\rm K}_1(\mS_n) \simeq K^*$ (Theorem
\ref{17Apr10}) is to use the fact that the group $\GL_\infty
(\mS_{n-1})$ is canonically isomorphic to the congruence subgroup
$(1+\gp_n)^*$ of $\mS_n^*=K^*\times (1+\ga_n)^*$,
$(1+\gp_n)^*\subseteq (1+\ga_n)^*$, where $\gp_n$ is an
(arbitrary) height one  prime  ideal of the algebra $\mS_n$. The
group $\mS_n^*$ is huge, {\em e.g.} 
\begin{equation}\label{Sa2nG}
 \mS_n^* \supset (1+\ga_n)^* \supset  \underbrace{\GL_\infty
(K)\ltimes\cdots \ltimes \GL_\infty (K)}_{2^n-1 \;\; {\rm times}},
\end{equation}
the iterated semi-direct product being a small part of the group
$\mS_n^*$.  The key ingredients in finding the groups $G_n$, $\Inn
(\mS_n)$ and $\mS_n^*$ (and their explicit generators)  are the
Fredholm operators and their indices, the current subgroups,  and
the ${\rm K}_1$-theory. This explains why it is possible to
recover the group $\GL_\infty (\mS_{n-1})$ in $\mS_n^*$ (this is
not straightforward),   to find its explicit generators. We prove in Theorem \ref{16Apr10}, Lemma \ref{a16Apr10}, and (\ref{cn1det})) that

\begin{itemize}
\item  {\em the group $\GL_\infty (\mS_n)$ is generated
by $K^*$, the group of elementary matrices $E_\infty (\mS_n)$ and
 $ (n-2)2^{n-1}+1$ matrices}  $\Bigl(
\begin{matrix} \th_{ij}(J)  &0 \\0 &1
\end{matrix}\Bigr)$
(Lemma \ref{a17Apr10}) where (see (\ref{ththdef}))
$$\th_{ij}(J):= \Bigl( 1+(y_i-1)\prod_{k\in J\backslash
i}(1-x_ky_k)\Bigr) \Bigl( 1+(x_j-1)\prod_{l\in J\backslash
j}(1-x_ly_l)\Bigr)\in (1+\ga_n)^*, $$ $J$ is a subset of $\{ 1,
\ldots , n\}$ with $|J|\geq 2$,  $i$ is the largest  number in $J$
and $j\in J\backslash i$.
\end{itemize} The final and the most difficult  part of the proof is to show
that
\begin{itemize}
\item (Theorem \ref{8May10}) {\em all the above matrices} $\Bigl(
\begin{matrix} \th_{ij}(J)  &0 \\0 &1
\end{matrix}\Bigr)$ {\em are elementary, i.e. the units $\th_{ij}(J)$ are elementary when regarded as matrices via the inclusion $\GL_1(\mS_n)\subseteq \GL_\infty (\mS_n)$}.  $\Box$
\end{itemize}
 We spend all of  Section \ref{THM8May10} to prove this fact.
\begin{itemize}
\item (Theorem \ref{13May10}) {\em Let $\gp$ be a nonzero
idempotent prime ideal of the algebra $\mS_n$ and $m={\rm ht} (\gp )$ be
its height. Then}
$$ {\rm K}_1(\mS_{n}, \gp )\simeq \begin{cases}
K^*, & \text{if }m=1,\\
\Z^{m\choose 2}\times K^{*m}& \text{if }m> 1.\\
\end{cases}$$
\end{itemize}

  Let $\Theta_{n,s}$, $s=1,
\ldots , n-1$, denote the  finitely generated subgroup of the group
$(1+\ga_n)^*$,  generated by the elements $\th_{ij}(J)$ where $J$ is
a subset of $\{ 1, \ldots , n\}$ with $|J|=s+1\geq 2$, and $i$ and
$j$ are two distinct elements of the set $J$. These, the so-called {\em current subgroups},  were introduced in \cite{K1aut} and \cite{Snaut}, and
they are the core (the non-obvious part) of the groups $G_n$,
$\Inn (\mS_n)$ and $\mS_n^*$ and the key for determining  the groups
$\GL_\infty (\mS_n)$, ${\rm K}_1(\mS_n)$, $\GL_\infty (\mS_n, \gp
)$ and  ${\rm K}_1(\mS_n, \gp)$,  as this paper demonstrates.

The paper is organized as follows. In Section \ref{KTG1S}, some
necessary results and constructions are collected for the algebra
$\mS_n$ and the group $(1+\ga_n)^*$. In Section \ref{EEE}, the
groups ${\rm K}_1(\mS_n)$, $\GL_\infty (\mS_n)$ and their explicit
generators are found. In Section \ref{THM8May10}, Theorem
\ref{8May10} is proved. In Section \ref{K1GPSN}, the groups
$\GL_\infty (\mS_n, \gp )$, ${\rm K}_1(\mS_n, \gp)$ and
explicit generators for them are found, and Theorem \ref{13May10} is
proved.

$\noindent $

{\bf The structure of the proof of Theorem \ref{13May10}}. The line of
proof of Theorem \ref{13May10} follows that  of
Theorem \ref{17Apr10} (but there are new moments): first, we prove
that the group $\GL_\infty (\mS_n, \gp )$ is generated by the
group $E_\infty (\mS_n , \gp )$ of $\gp$-elementary matrices, some
explicit `diagonal' matrices, and {\em some} of the matrices
 $\Bigl(\begin{matrix} \th_{ij}(J)  &0 \\0 &1
\end{matrix}\Bigr)$  (Theorem \ref{12May10}, Lemma
\ref{b13May10}). Then an `obvious' normal subgroup $\CE (\mS_n ,
\gp )$ of $\GL_\infty ( \mS_n , \gp )$ is introduced and we prove
that
$$ \GL_\infty (\mS_n , \gp )/\CE (\mS_n ,
\gp )\simeq\begin{cases}
K^*, & \text{if }m=1,\\
\Z^{m\choose 2}\times K^{*m}& \text{if }m> 1.\\
\end{cases}$$
This gives  the inclusion $E_\infty (\mS_n, \gp ) \subseteq \CE
(\mS_n , \gp )$. The key moment in proving that the  opposite
inclusion holds  is (surprisingly) the fact that ${\rm K}_1(\mS_n)
\simeq K^*$. The new moment is that not all the `diagonal'
matrices and not all the matrices $\Bigl(\begin{matrix}
\th_{ij}(J)  &0 \\0 &1
\end{matrix}\Bigr)$  that form a part of the generating set for
the group $\GL_\infty (\mS_n , \gp )$ are $\gp $-elementary. $\Box
$

A canonical form is found (Theorem \ref{13May10}) for each element
$a\in \GL_\infty (\mS_n , \gp )$. Using it,  an effective criterion
(Corollary \ref{b15May10}) is given for an element $a\in
\GL_\infty (\mS_n , \gp )$ to be a product of $\gp$-elementary
matrices, i.e. $a\in E_\infty (\mS_n , \gp )$.


\section{The groups $\mS_n^*$ and  $(1+\ga_n)^*$ and their subgroups}\label{KTG1S}

In this section, we collect some results without proofs on the
algebras $\mS_n$ from \cite{shrekalg} and \cite{Snaut} that will
be used in this paper, their proofs can be found in
\cite{shrekalg} and  \cite{Snaut}.
 Several important
subgroups of the group $(1+\ga_n)^*$ are considered. The most
interesting of these are the current subgroups $\Theta_{n,s}$,
$s=1, \ldots , n-1$. They encapsulate the most difficult parts of
the groups $\mS_n^*$ and $G_n$. 

$\noindent $

{\bf The algebra of one-sided inverses of a polynomial algebra}.
Clearly, $\mathbb{S}_n=\mS_1(1)\t \cdots \t\mS_1(n)\simeq
\mathbb{S}_1^{\t n}$ where $\mS_1(i):=K\langle x_i,y_i \, | \,
y_ix_i=1\rangle \simeq \mS_1$ and $\mS_n=\bigoplus_{\alpha , \beta
\in \N^n} Kx^\alpha y^\beta$ where $x^\alpha := x_1^{\alpha_1}
\cdots x_n^{\alpha_n}$, $\alpha = (\alpha_1, \ldots , \alpha_n)$,
$y^\beta := y_1^{\beta_1} \cdots y_n^{\beta_n}$ and $\beta =
(\beta_1,\ldots , \beta_n)$. In particular, the algebra $\mS_n$
contains two polynomial subalgebras $P_n$ and $Q_n:=K[y_1, \ldots
, y_n]$ and is equal,  as a vector space,  to their tensor product
$P_n\t Q_n$.

When $n=1$, we usually drop the subscript `1' if this does not
lead to confusion.  So, $\mS_1= K\langle x,y\, | \,
yx=1\rangle=\bigoplus_{i,j\geq 0}Kx^iy^j$. For each natural number
$d\geq 1$, let $M_d(K):=\bigoplus_{i,j=0}^{d-1}KE_{ij}$ be the
algebra of $d$-dimensional matrices where $\{ E_{ij}\}$ are the
matrix units, and
$$M_\infty (K) :=
\varinjlim M_d(K)=\bigoplus_{i,j\in \N}KE_{ij}$$ is the algebra
(without 1) of infinite dimensional matrices. The algebra $\mS_1$
contains the ideal $F:=\bigoplus_{i,j\in \N}KE_{ij}$, where
\begin{equation}\label{Eijc}
E_{ij}:= x^iy^j-x^{i+1}y^{j+1}, \;\; i,j\geq 0.
\end{equation}
For all natural numbers $i$, $j$, $k$, and $l$,
$E_{ij}E_{kl}=\d_{jk}E_{il}$ where $\d_{jk}$ is the Kronecker
delta function.  The ideal $F$ is an algebra (without 1)
isomorphic to the algebra (without 1)  $M_\infty (K)$ via
$E_{ij}\mapsto E_{ij}$.  For all $i,j\geq 0$, 
\begin{equation}\label{xyEij}
xE_{ij}=E_{i+1, j}, \;\; yE_{ij} = E_{i-1, j}\;\;\; (E_{-1,j}:=0),
\end{equation}
\begin{equation}\label{xyEij1}
E_{ij}x=E_{i, j-1}, \;\; E_{ij}y = E_{i, j+1} \;\;\;
(E_{i,-1}:=0).
\end{equation}
The algebra
\begin{equation}\label{mS1d}
\mS_1= K\oplus xK[x]\oplus yK[y]\oplus F
\end{equation}
is a direct sum of vector spaces. Then 
\begin{equation}\label{mS1d1}
\mS_1/F\simeq K[x,x^{-1}]=:L_1, \;\; x\mapsto x, \;\; y \mapsto
x^{-1},
\end{equation}
since $yx=1$, $xy=1-E_{00}$ and $E_{00}\in F$.

$\noindent $

The algebra $\mS_n = \bigotimes_{i=1}^n \mS_1(i)$ contains the
ideal
$$F_n:= F^{\t n }=\bigoplus_{\alpha , \beta \in
\N^n}KE_{\alpha \beta}, \;\; {\rm where}\;\; E_{\alpha
\beta}:=\prod_{i=1}^n E_{\alpha_i \beta_i}(i), \;\; E_{\alpha_i
\beta_i}(i):=x_i^{\alpha_i}y_i^{\beta_i}-x_i^{\alpha_i+1}y_i^{\beta_i+1}.$$
Note that $E_{\alpha \beta}E_{\g \rho}=\d_{\beta \g }E_{\alpha
\rho}$ for all elements $\alpha, \beta , \g , \rho \in \N^n$ where
$\d_{\beta
 \g }$ is the Kronecker delta function;  $F_n=\bigotimes_{i=1}^nF(i)$ and
 $F(i):=\bigoplus_{s,t\in \N}KE_{st}(i)$.

\begin{itemize}
 \item   {\em The
 algebra $\mS_n$ is central, prime and catenary. Every nonzero
 ideal  of $\mS_n$ is an essential left and right submodule of}
 $\mS_n$.
 \item  {\em The ideals of
 $\mS_n$ commute ($IJ=JI$);  and the set of ideals of $\mS_n$ satisfy the a.c.c..}
\item {\em $\ga \cap \gb = \ga\gb$ for all idempotent ideals $\ga
$ and $\gb$ of the algebra} $\mS_n$; \item  {\em The classical
Krull dimension $\clKdim (\mS_n)$ of $\mS_n$ is $2n$.}
  \item  {\em Let $I$
be an ideal of $\mS_n$. Then the factor algebra $\mS_n / I$ is
left (or right) Noetherian iff the ideal $I$ contains all the
height one prime ideals  of the algebra $\mS_n$.}
\end{itemize}

{\bf The set of height one prime ideals of $\mS_n$}.   Consider
the ideals of the algebra $\mS_n$:
$$\gp_1:=F\t \mS_{n-1},\; \gp_2:= \mS_1\t F\t \mS_{n-2}, \ldots ,
 \gp_n:= \mS_{n-1} \t F.$$ Then $\mS_n/\gp_i\simeq
\mS_{n-1}\t (\mS_1/F) \simeq  \mS_{n-1}\t K[x_i, x_i^{-1}]$ and
$\bigcap_{i=1}^n \gp_i = \prod_{i=1}^n \gp_i =F^{\t n }=F_n$.
Clearly, $\gp_i\not\subseteq \gp_j$ for all $i\neq j$.

\begin{itemize}
 \item
{\em The set $\CH_1$ of height one prime ideals of the algebra
$\mS_n$ is} $\{ \gp_1, \ldots , \gp_n\}$.
\end{itemize}

 Let
$\ga_n:= \gp_1+\cdots +\gp_n$. Then the factor algebra
\begin{equation}\label{SnSn}
\mS_n/ \ga_n\simeq (\mS_1/F)^{\t n } \simeq \bigotimes_{i=1}^n
K[x_i, x_i^{-1}]= K[x_1, x_1^{-1}, \ldots , x_n, x_n^{-1}]=:L_n
\end{equation}
is a skew Laurent polynomial algebra in $n$ variables,  and so
$\ga_n$ is a prime ideal of height and co-height $n$ of the
algebra $\mS_n$.

\begin{proposition}\label{a19Dec8}
{\rm \cite{shrekalg}} The polynomial algebra $P_n$
 is the only (up to isomorphism)  faithful  simple $\mS_n$-module.
\end{proposition}

In more detail, ${}_{\mS_n}P_n\simeq \mS_n / (\sum_{i=0}^n \mS_n
y_i) =\bigoplus_{\alpha \in \N^n} Kx^\alpha \overline{1}$,
$\overline{1}:= 1+\sum_{i=1}^n \mS_ny_i$; and the action of the
canonical generators of the algebra $\mS_n$ on the polynomial
algebra $P_n$ is given by the rule:
$$ x_i*x^\alpha = x^{\alpha + e_i}, \;\; y_i*x^\alpha = \begin{cases}
x^{\alpha - e_i}& \text{if } \; \alpha_i>0,\\
0& \text{if }\; \alpha_i=0,\\
\end{cases}  \;\; {\rm and }\;\; E_{\beta \g}*x^\alpha = \d_{\g
\alpha} x^\beta,
$$
where the set $e_1:= (1,0,\ldots , 0),  \ldots , e_n:=(0, \ldots ,
0,1)$ is the canonical basis for the free $\Z$-module
$\Z^n=\bigoplus_{i=1}^n \Z e_i$.  We identify the algebra $\mS_n$
with its image in the algebra $\End_K(P_n)$ of all the $K$-linear
maps from the vector space $P_n$ to itself, i.e. $\mS_n \subset
\End_K(P_n)$.

For each non-empty subset $I$ of the set $\{ 1,\ldots , n\}$, let
$\mS_I:=\bigotimes_{i\in I}\mS_1(i)\simeq \mS_{|I|}$  where $|I|$
is the number of elements in the set $I$, $F_I:=\bigotimes_{i\in
I}F(i)\simeq M_\infty (K)$, $\ga_I$ is the ideal of the algebra
$\mS_I$ generated by the vector space $\bigoplus_{i\in I}F(i)$,
i.e. $\ga_I:=\sum_{i\in I}F(i)\t \mS_{I\backslash i}$. The factor
algebra $L_I:=\mS_I/\ga_I\simeq K[x_i,x_i^{-1}]_{i\in I}$ is a
Laurent polynomial algebra. For elements $\alpha =
(\alpha_i)_{i\in I}, \beta = (\beta_i)_{i\in I}\in \N^I$, let $
E_{\alpha \beta}(I):=\prod_{i\in I}E_{\alpha_i\beta_i}(i)$. Then
$E_{\alpha\beta}(I)E_{\xi\rho}(I)= \d_{\beta
\xi}E_{\alpha\rho}(I)$ for all $\alpha , \beta , \xi , \rho \in
\N^I$.

$\noindent $

{\bf The $G_n$-invariant normal subgroups $(1+\ga_{n,s})^*$  of
$(1+\ga_n)^*$.} Let $G_n:=\Aut_{K-{\rm alg}}(\mS_n)$.  We will use
often the following obvious lemma.

\begin{lemma}\label{aa13Dec8}
{\rm \cite{shrekaut}} Let $R$ be a ring and $I_1, \ldots , I_n$ be
ideals of the ring $R$ such that $I_iI_j=0$ for all $i\neq j$. Let
$a= 1+a_1+\cdots +a_n\in R$ where $a_1\in I_1, \ldots , a_n\in
I_n$. The element $a$ is a unit of the ring $R$ iff all the
elements $1+a_i$ are units; and, in this case,  $a^{-1}=
(1+a_i)^{-1} (1+a_2)^{-1}\cdots (1+a_n)^{-1}$.
\end{lemma}

Let $R$ be a ring, $R^*$ be its group of units, $I$ be an ideal of
$R$ such that $I\neq R$, and let  $ (1+I)^*$ be the group of units
of the multiplicative monoid $1+I$. Then $R^*\cap (1+I)= (1+I)^*$
 and $(1+I)^*$ is a normal subgroup of $R^*$.

For each subset $I$ of the set $\{ 1, \ldots , n\}$, let $\gp_I:=
\bigcap_{i\in I}\gp_i$, and $\gp_\emptyset :=\mS_n$. Each $\gp_I$
is an ideal of the algebra $\mS_n$ and $\gp_I=\prod_{i\in
I}\gp_i$.  The complement to the subset $I$ is denoted by $CI$.
For a one-element subset $\{ i\}$, we write $Ci$ rather than $C\{
i\}$. In particular, $\gp_{Ci}:= \gp_{C\{ i\} }=\bigcap_{j\neq
i}\gp_j$.

 For each number $s=1,\ldots , n$, let
$\ga_{n,s}:=\sum_{|I|=s}\gp_I$. By the very definition, the ideals
$\ga_{n,s}$ are $G_n$-invariant ideals (since the set $\CH_1$ of
all the height one prime ideals of the algebra $\mS_n$ is $\{
\gp_1, \ldots , \gp_n\}$, \cite{shrekaut}, and $\CH_1$ is a
$G_n$-orbit). We have a strictly descending chain  of
$G_n$-invariant  ideals of the algebra $\mS_n$:
$$ \ga_n=\ga_{n,1}\supset \ga_{n,2}\supset \cdots \supset
\ga_{n,s}\supset \cdots \supset \ga_{n,n}=F_n\supset
\ga_{n,n+1}:=0.$$ These are also ideals of the subalgebra
$K+\ga_n$ of $\mS_n$. Each set $\ga_{n,s}$ is an ideal of the
algebra $K+\ga_{n,t}$ for all $t\leq s$, and the group of units of
the algebra $K+\ga_{n,s}$ is the direct product of its two
subgroups
$$ (K+\ga_{n,s})^* = K^*\times (1+\ga_{n,s})^*, \;\; s=1, \ldots ,
n.$$ The groups $(K+\ga_{n,s})^*$ and $(1+\ga_{n,s})^*$ are
$G_n$-invariant. There is the descending chain of $G_n$-invariant
(hence normal) subgroups of $( 1+\ga_n)^*$:
$$ (1+\ga_n)^*=(1+\ga_{n,1})^*\supset \cdots \supset
(1+\ga_{n,s})^*\supset \cdots \supset
(1+\ga_{n,n})^*=(1+F_n)^*\supset (1+\ga_{n,n+1})^*=\{ 1\}.$$
 For each number $s=1, \ldots , n$, the factor
algebra
$$(K+\ga_{n,s})/\ga_{n,s+1}=K\bigoplus
\bigoplus_{|I|=s}\bgp_I$$ contains the idempotent ideals
$\bgp_I:=(\gp_I+\ga_{n,s+1})/\ga_{n,s+1}$ such that $\bgp_I
\bgp_J=0$ for all $I\neq J$ such that $|I|=|J|=s$.

Recall that for a Laurent polynomial algebra $L=K[x_1^{\pm 1},
\ldots , x_n^{\pm 1}]$, ${\rm K}_1(L)\simeq L^*$, \cite{Swan},
\cite{Bass-book-K-theory}, \cite{Milnor-book-K-theory},
\begin{equation}\label{GLLU}
\GL_\infty (L)=U(L)\ltimes E_\infty (L)
\end{equation}
where $E_\infty (L)$ is the subgroup of $\GL_\infty (L)$ generated
by all the {\em elementary matrices} $\{ 1+aE_{ij}\, | \, a\in L,
i,j\in \N, i\neq j\}$, and $U(L):=\{\mu (u):= uE_{00}+1-E_{00}\, |
\, u\in L^*\}\simeq L^*$, $\mu (u)\lra u$. The group $E_\infty
(L)$ is a normal subgroup of $\GL_\infty (L)$. This is true for an
arbitrary coefficient ring.

By Lemma \ref{aa13Dec8} and (\ref{GLLU}), the group of units of
the algebra $(K+\ga_{n,s})/\ga_{n,s+1}=:K+\ga_{n,s}/\ga_{n,s+1}$
is the direct product of groups,
$$ (K+\ga_{n,s}/\ga_{n,s+1})^*=K^*\times
\prod_{|I|=s}(1+\bgp_I)^*\simeq K^*\times \prod_{|I|=s}\GL_\infty
(L_{CI})\simeq K^*\times \prod_{|I|=s}U(L_{CI})\ltimes E_\infty
(L_{CI})$$ since $(1+\bgp_I)^*\simeq (1+M_\infty
(L_{CI}))^*=\GL_\infty (L_{CI})$ where
$L_{CI}:=\mS_{CI}/\ga_{CI}=\bigotimes_{i\in CI}K[x_i, x_i^{-1}]$
is the Laurent polynomial algebra.
 In more detail, for each non-empty subset $I$ of $\{ 1, \ldots ,
 n\}$, let $\Z^I:=\bigoplus_{i\in I}\Z e_i$. It is a subgroup of
 $\Z^n=\bigoplus_{i=1}^n\Z e_i$. Similarly, $\N^I:= \bigoplus_{i\in I}\N e_i$. By (\ref{GLLU}),
\begin{equation}\label{apns}
(1+\bgp_I)^*=U(L_{CI})\ltimes E_\infty (L_{CI})=(U_I(K)\times
\mX_{CI})\ltimes E_\infty (L_{CI})
\end{equation}
where
\begin{eqnarray*}
 U(L_{CI})&:=&\{\mu_I(u):= uE_{00}(I)+1-E_{00}(I)\, | \, u\in
L^*_{CI}\}\simeq L_{CI}^*,\; \mu_I(u)\lra u, \\
L^*_{CI}&=&\{\l x^\alpha \, | \, \l \in K^*, \alpha \in \Z^{CI}\},\\
U_I(K)&:=&\{ \mu_I(\l ):=\l E_{00}(I)+1-E_{00}(I)\, |\, \l \in
K^*\} \simeq
K^*,\;  \mu_I(\l )\lra \l, \\
\mX_{CI}&:=&\{ \mu_I(x^\alpha ):=x^\alpha E_{00}(I)+1-E_{00}(I)\,
| \, \alpha \in
\Z^{CI}\}\simeq \Z^{CI}\simeq \Z^{n-s},\;  \mu_I(x^\alpha )\lra \alpha , \\
E_\infty (L_{CI})&:=&\langle 1+aE_{\alpha\beta}(I)\, | \, a\in
L_{CI}, \alpha,\beta \in \N^I, \alpha\neq \beta \rangle.
\end{eqnarray*}

The algebra  epimorphism  $\psi_{n,s}: K+\ga_{n,s}\ra
(K+\ga_{n,s})/\ga_{n,s+1}$, $ a\mapsto a+\ga_{n,s+1}$, yields the
group homomorphism of their groups of units $(K+\ga_{n,s})^*\ra
(K+\ga_{n,s}/\ga_{n,s+1})^*$ and whose kernel is
$(1+\ga_{n,s+1})^*$. As a result we have an  exact sequence of
group homomorphisms: 
\begin{equation}\label{apns1}
1\ra (1+\ga_{n,s+1})^*\ra (1+\ga_{n,s})^* \stackrel{\psi_{n,s}}\ra
\prod_{|I|=s}(1+\bgp_I)^* \simeq \prod_{|I|=s}\GL_\infty
(L_{CI})\ra \CZ_{n,s}\ra 1.
\end{equation}
For $s=n$, the map $\psi_{n,n}$ is the identity map, and so
$\CZ_{n,n}=\{ 1\}$.  Intuitively, the group $\CZ_{n,s}$ represents
`relations' that determine the image $\im (\psi_{n,s})$ as a
subgroup of
 $\prod_{|I|=s}(1+\bgp_I)^*$. The group
$\CZ_{n,s}$ is a free abelian group of rank ${n\choose s+1}$,
\cite{Snaut}. So, the image of the map $\psi_{n,s}$ is large. Note
that $\ga_{n,s+1}$ and $\gp_I$ (where $|I|=s)$ are ideals of the
algebra $K+\ga_{n,s}$. The groups $(1+\ga_{n,s+1})^*$ and
$(1+\gp_I)^*$ (where $|I|=s$) are normal subgroups of
$(1+\ga_{n,s})^*$. Thus  the  subgroup $\Upsilon_{n,s}$ of
$(1+\ga_{n,s})^*$ generated by these normal subgroups is a normal
subgroup of $(1+\ga_{n,s})^*$. As a subset of $(1+\ga_{n,s})^*$,
the group $\Upsilon_{n,s}$ is equal to the product of the groups
$(1+\ga_{n,s+1})^*$, $(1+\gp_I)^*$, $|I|=s$, in {\em arbitrary}
order (by their normality), i.e. 
\begin{equation}\label{UPns}
\Upsilon_{n,s}=\prod_{|I|=s}(1+\gp_I)^*\cdot (1+\ga_{n,s+1})^*.
\end{equation}
By Theorem \ref{I1}, the group $\Upsilon_{n,s}$ is a
$G_n$-invariant (hence, normal) subgroup of $\mS_n^*$.  The factor
group $(1+\ga_{n,s})^*/ \Upsilon_{n,s}$ is a free abelian group of
rank ${n\choose s+1}s$, \cite{Snaut}.

By (\ref{apns}), the direct product of groups
$\prod_{|I|=s}(1+\bgp_I)^*=\mX_{n,s}\ltimes \bG_{n,s}$ is the
semi-direct product of its two subgroups 
\begin{equation}\label{apns2}
\mX_{n,s}:=\prod_{|I|=s}\mX_{CI}\simeq \Z^{{n\choose s}(n-s)}\;\;
{\rm and}\;\;  \bG_{n,s}:=\prod_{|I|=s}U_I(K)\ltimes
E_\infty(L_{CI}).
\end{equation}
For each subset $I$ of $\{ 1,\ldots , n\}$ such that $|I|=s$,
$U_I(K)\ltimes E_\infty (\mS_{CI})$ is a subgroup of
$(1+\gp_I)^*$ where 
\begin{equation}\label{apns3}
U_I(K):=\{ \mu_I(\l )\, | \, \l \in K^*\} \simeq K^*, \; E_\infty
(\mS_{CI}):=\langle 1+aE_{\alpha\beta}(I)\, | \, a\in \mS_{CI},
\alpha \neq  \beta \in \N^I\rangle ,
\end{equation}
where $\mu_I(\l ):=\l E_{00}(I)+1-E_{00}(I)$.  Clearly,
$$ \psi_{n,s}|_{U_I(K)}: U_I(K)\simeq U_I(K),\;\;
\mu_I(\l )\mapsto \mu_I(\l ), $$ and $\psi_{n,s}(U_I(K)\ltimes
E_\infty (\mS_{CI}))=U_I(K)\ltimes E_\infty (L_{CI})$ for all
subsets $I$ with $|I|=s$. The subgroup of $(1+\ga_{n,s})^*$,
\begin{equation}\label{apns4}
\G_{n,s}:=\psi^{-1}_{n,s}(\bG_{n,s})={}^{set}\prod_{|I|=s}(U_I(K)\ltimes
E_\infty (\mS_{CI}))\cdot (1+\ga_{n,s+1})^*,
\end{equation}
is a normal subgroup as it is  the  pre-image of a normal subgroup. We
added the upper script `set' to indicate that this is a product of
subgroups but not a direct product, in general. It is obvious
that $\psi_{n,s}(\G_{n,s})=\bG_{n,s}$ and $\G_{n,s}\subseteq
\Upsilon_{n,s}$. In fact, $\G_{n,s}= \Upsilon_{n,s}$,
\cite{Snaut}.
  Let
$\D_{n,s}:=(1+\ga_{n,s})^*/\G_{n,s}$. The group homomorphism
$\psi_{n,s}$ (see (\ref{apns1}))  induces the group monomorphism
$$ \bpsi_{n,s}:\D_{n,s}\ra
\prod_{|I|=s}(1+\bgp_I)^*)/\bG_{n,s}\simeq \mX_{n,s}\simeq
\Z^{{n\choose s}(n-s)}.$$ This means that the group $\D_{n,s}$ is
a free abelian group of rank $\leq {n\choose s}(n-s)$. In fact,
the rank is equal to ${n\choose s+1}s$, \cite{Snaut}.

For each subset $I$ with $|I|=s$, consider the  free abelian group
$\mX_{CI}':= \bigoplus_{j\in CI}\Z (j,I)\simeq \Z^{n-s}$ where $\{
(j,I)\, | \, j\in CI\}$ is its free basis. Let
$$ \mX_{n,s}':=
\bigoplus_{|I|=s}\mX_{CI}'=\bigoplus_{|I|=s}\bigoplus_{j\in CI}\Z
(j,I)\simeq \Z^{{n\choose s}(n-s)}.$$ For each subset $I$,
consider the isomorphism of abelian groups
$$ \mX_{CI}\ra \mX_{CI}', \;\; \mu_I(x_j):= x_jE_{00}(I)+1-E_{00}(I)\mapsto (j,I).
$$
These isomorphisms yield the group isomorphism 
\begin{equation}\label{XXns}
\mX_{n,s}\ra \mX_{n,s}', \;\; \mu_I(x_j)\mapsto (j,I).
\end{equation}
Each element $a$ of the group $\mX_{n,s}$ is a unique product
$a=\prod_{|I|=s}\prod_{j\in CI}\mu_I(x_j)^{n(j,I)}$ where
$n(j,I)\in \Z$. Each element $a'$ of the group $\mX_{n,s}'$ is a
unique sum $a'=\sum_{|I|=s}\sum_{j\in CI}n(j,I)\cdot (j,I)$ where
$n(j,I)\in \Z$. The map (\ref{XXns}) sends $a$ to $a'$. To make
computations more readable we set $e_I:=E_{00}(I)$. Then
$e_Ie_J=e_{I\cup J}$.

$\noindent $

{\bf The current groups $\Theta_{n,s}$, $s=1,\ldots , n-1$}. The
current groups $\Theta_{n,s}$ are the most important subgroups of
the group $(1+\ga_n)^*$. They are finitely generated groups and
generators are given explicitly.  The generators of the groups
$\Theta_{n,s}$ are units of the algebra $\mS_n$ but they are
defined as a product of two {\em non-units}. As a result the
groups $\Theta_{n,s}$ capture the most delicate phenomena regarding
the structure and  properties of the groups $\mS_n^*$ and
$G_n$.

For each non-empty subset $I$ of $\{ 1,\ldots , n\}$ with
$s:=|I|<n$ and an element $i\in CI$, let
$$X(i,I):=\mu_I(x_i)= x_iE_{00}(I)+1-E_{00}(I)\;\; {\rm and}\;\;
Y(i,I):=\mu_I(y_i)= y_iE_{00}(I)+1-E_{00}(I).$$ Then
$Y(i,I)X(i,I)=1$, $\ker \, Y(i,I)=P_{C(I\cup i)}$, and $P_n=\im \,
X(i,I)\bigoplus P_{C(I\cup i)}$ where $P_{C(I\cup i)}:=
K[x_j]_{j\in C(I\cup i)}$. Recall that $\mS_n\subset \End_K(P_n)$.
As an element of the algebra $\End_K(P_n)$, the map $X(i,I)$ is
injective (but not bijective), and the map $Y(i,I)$ is surjective
(but not bijective).

$\noindent $

{\it Definition}. For each subset $J$ of $\{ 1,\ldots , n\}$ with
$|J|=s+1\geq 2$ and for two distinct elements $i$ and $j$ of the
set $J$, let
\begin{equation}\label{ththdef}
 \th_{ij}(J):=Y(i,J\backslash i) X(j,J\backslash j)\in
(1+\gp_{J\backslash i}+\gp_{J\backslash j})^* \subseteq
(1+\ga_{n,s})^*.
\end{equation}
The {\bf current group} $\Theta_{n,s}$ is the
subgroup of $(1+\ga_{n,s})^*$ generated by all the elements
$\th_{ij}(J)$ (for all the possible choices of $J$, $i,$ and $j$).

$\noindent $

The unit $\th_{ij}(I)$ is the product in $\End_K(P_n)$ of an
injective map and a surjective map,  none of which is a  bijection.
\begin{equation}\label{tiji}
\th_{ij}(J)= \th_{ji}(J)^{-1}.
\end{equation}

Suppose that $i$, $j$, and $k$ are distinct elements of the set
$J$ (hence $|J|\geq 3$). Then 
\begin{equation}\label{tijjk}
\th_{ij}(J)\th_{jk}(J)=\th_{ik}(J).
\end{equation}

For each number $s=1, \ldots , n-1$, the free abelian group
$\mX_{n,s}'$ admits the decomposition
$\mX_{n,s}'=\bigoplus_{|J|=s+1}\bigoplus_{j\cup I=J}\Z (j,I)$, and
using it we define a character (a homomorphism) $\chi_J'$, for
each subset $J$ with $|J|=s+1$:
$$ \chi_J':\mX_{n,s}'\ra \Z, \;\; \sum_{|J'|=s+1}\sum_{j\cup I=J'}n_{j,I}(j,I)\mapsto \sum_{j\cup I=J}n_{j,I}.$$
Let $\max (J)$ be the largest number  in the set $J$. The group
$\mX_{n,s}'$ is the direct sum 
\begin{equation}\label{XK1}
\mX_{n,s}'=\mK_{n,s}'\bigoplus \mY_{n,s}'
\end{equation}
of its free abelian subgroups,
\begin{eqnarray*}
 \mK_{n,s}'&=& \bigcap_{|J|=s+1}\ker (\chi_J')=\bigoplus_{|J|=s+1}\bigoplus_{j\in J\backslash
 \max (J)}\Z (-(\max (J), J\backslash \max (J))+(j,J\backslash j))\simeq \Z^{{n\choose s+1}s},  \\
 \mY_{n,s}'&=& \bigoplus_{|J|=s+1}\Z (\max (J), J\backslash \max
 (J))\simeq \Z^{{n\choose s+1}}.
\end{eqnarray*}
The same decompositions hold if instead of $\max (J)$,  we choose
any element of the set $J$.  Consider the group homomorphism
$\psi_{n,s}':(1+\ga_{n,s})^*\ra \mX_{n,s}'$ defined as the
composition of the following group homomorphisms:
$$\psi_{n,s}': (1+\ga_{n,s})^*\ra
(1+\ga_{n,s})^*/\G_{n,s}\stackrel{\bpsi_{n,s}}{\ra}
\prod_{|I|=s}(1+\bgp_I)^*/\bG_{n,s}\simeq \mX_{n,s}\simeq
\mX_{n,s}'.$$ Then 
\begin{equation}\label{ptijJ}
\psi_{n,s}'(\th_{ij}(J))=-(i,J\backslash i)+(j,J\backslash j).
\end{equation}
It follows that 
\begin{equation}\label{ptijJ1}
\psi_{n,s}'(\Theta_{n,s})=\mK_{n,s}',
\end{equation}
since, by (\ref{ptijJ}), $\psi_{n,s}'(\Theta_{n,s})\supseteq
\mK_{n,s}'$ (as the free basis for $\mK_{n,s}'$, introduced above,
belongs to the set $\psi_{n,s}'(\Theta_{n,s})$); again, by
(\ref{ptijJ}), $\psi_{n,s}'(\Theta_{n,s})\subseteq
\bigcap_{|J|=s+1}\ker (\chi_J')= \mK_{n,s}'$.

Let $H, H_1, \ldots , H_m$ be subsets (usually subgroups) of a
group $H$.  We say that $H$ is the {\em  product} of $H_1, \ldots
, H_m$, and write $H={}^{set}\prod_{i=1}^mH_i=H_1\cdots H_m$, if
each element $h$ of $H$ is a  product $h=h_1\cdots h_m$ where
$h_i\in H_i$. We add the subscript `set' (sometime) in order to
distinguish it from the direct product of groups.
 We say that $H$ is the {\em exact product} of $H_1,
\ldots , H_m$, and write
$H={}^{exact}\prod_{i=1}^mH_i=H_1\times_{ex}\cdots
\times_{ex}H_m$, if each element $h$ of $H$ is a {\em unique}
product $h=h_1\cdots h_m$ where $h_i\in H_i$. The order in the
definition of the exact product is important. A {\em semi-direct product} of groups $H_1, \ldots , H_m$ is denoted by
$$ H_1\ltimes (H_2\ltimes (\cdots\ltimes H_m))= H_1\ltimes H_2\ltimes \cdots \ltimes H_m = {}^{semi} \prod_{i=1}^m H_i. $$
The subgroup of $(1+\ga_{n,s})^*$ generated by the groups
$\Theta_{n,s}$ and $\G_{n,s}$ is equal to their product
$\Theta_{n,s}\G_{n,s}$, by the normality of $\G_{n,s}$. The
subgroup $\G_{n,s}$ of the group $\Theta_{n,s}\G_{n,s}$ is a
normal subgroup, hence the intersection $\Theta_{n,s}\cap
\G_{n,s}$ is a normal subgroup of $\Theta_{n,s}$.

\begin{lemma}\label{a24May9}
\cite{Snaut} For each number $s=1, \ldots , n-1$, the group
$\Theta_{n,s}\G_{n,s}$ is the semi-direct  product
$$\Theta_{n,s}\G_{n,s}= {}^{semi}\prod_{|J|=s+1}\prod_{j\in J\backslash \max (J)}\langle
\th_{\max (J), j}(J)\rangle\ltimes \G_{n,s},$$ where the order in
the double product is arbitrary. Each element $a\in
\Theta_{n,s}\G_{n,s}$ is a  unique product
$a=\prod_{|J|=s+1}\prod_{j\in J\backslash \max (J)} \th_{\max (J),
j}(J)^{n(j,J)}\cdot \g$ where $n(j,J)\in \Z$ and $\g \in
\G_{n,s}$.
\end{lemma}

$\noindent $

For each number $s=1, \ldots , n-1$, consider the subset of
$(1+\ga_{n,s})^*$, 
\begin{equation}\label{Tpns}
\Theta_{n,s}':= {}^{exact}\prod_{|J|=s+1}\prod_{j\in J\backslash
\max (J)}\langle \th_{\max (J), j}(J)\rangle,
\end{equation}
 which is the exact
product of cyclic groups (each of them is isomorphic to $\Z $),
since  each element $u$ of $\Theta_{n,s}'$ is a {\em unique}
product $u= \prod_{|J|=s+1}\prod_{j\in J\backslash \max
(J)}\th_{\max (J), j}(J)^{n(j,J)}$ where $n(j,J)\in \Z$ (Lemma
\ref{a24May9}).

By Lemma \ref{a24May9}, $\Theta_{n,s}/ \Theta_{n,s}\cap
\G_{n,s}\simeq \Theta_{n,s}\G_{n,s}/ \G_{n,s}\simeq
\mK_{n,s}'\simeq \Z^{{n\choose s+1}s}$, and so the commutant of
the current group $\Theta_{n,s}$ belongs  to the group $\G_{n,s}$,
i.e. 
\begin{equation}\label{TnsG}
[\Theta_{n,s}, \Theta_{n,s}]\subseteq \G_{n,s}.
\end{equation}
Recall that the {\em commutant} $[G,G]$ of a group $G$ is the
subgroup of $G$ generated by all  {\em group commutators}
$[a,b]:=aba^{-1}b^{-1}$ where $a,b\in G$. The commutant is a
normal subgroup.  The next theorem is the key point in finding
explicit generators for the groups $\mS_n^*$ and $G_n$.

\begin{theorem}\label{24May9}
\cite{Snaut}
$\psi_{n,s}'((1+\ga_{n,s})^*)=\psi_{n,s}'(\Theta_{n,s})$ for
$s=1,\ldots , n-1$.
\end{theorem}

For each number $s=1, \ldots , n-1$, consider the following
subsets of the group $(1+\ga_{n,s})^*$, 
\begin{equation}\label{mEns}
\mE_{n,s}:=\prod_{|I|=s}U_I(K)\ltimes E_\infty (\mS_{CI}) \;\;
{\rm and}\;\;  \mP_{n,s}:=\prod_{|I|=s}(1+\gp_i)^*.
\end{equation}
These are products of subgroups $(1+\ga_{n,s})^*$ in arbitrary order, but  which is
fixed for each $s$.

\begin{theorem}\label{A10May9}
\cite{Snaut}
\begin{enumerate}
\item $(1+\ga_n)^* = \Theta_{n,1}\G_{n,1} = \Theta_{n,1} \mE_{n,1}
\Theta_{n,2} \mE_{n,2}\cdots \Theta_{n,n-1} \mE_{n,n-1} $.
Moreover,  for $s=1, \ldots , n-1$,  $(1+\ga_{n,s})^* =
\Theta_{n,s}\G_{n,s} = \Theta_{n,s} \mE_{n,s} \Theta_{n,s+1}
\mE_{n,s+1}\cdots \Theta_{n,n-1} \mE_{n,n-1} $. \item $(1+\ga_n)^*
= \Theta_{n,1}\Upsilon_{n,1}=
 \Theta_{n,1} \mP_{n,1} \Theta_{n,2}
\mP_{n,2}\cdots \Theta_{n,n-1} \mP_{n,n-1} $. Moreover,  for $s=1,
\ldots , n-1$, $(1+\ga_{n,s})^* = \Theta_{n,s}\Upsilon_{n,s}=
 \Theta_{n,s} \mP_{n,s} \Theta_{n,s+1}
\mP_{n,s+1}\cdots \Theta_{n,n-1} \mP_{n,n-1} $.
\end{enumerate}
\end{theorem}

\begin{theorem}\label{B10May9}
\cite{Snaut}
\begin{enumerate}
\item $(1+\ga_n)^* = \Theta_{n,1}' \mE_{n,1} \Theta_{n,2}'
\mE_{n,2}\cdots \Theta_{n,n-1}' \mE_{n,n-1} $.  Moreover,  for
$s=1, \ldots , n-1$, $(1+\ga_{n,s})^* = \Theta_{n,s}' \mE_{n,s}
\Theta_{n,s+1}' \mE_{n,s+1}\cdots \Theta_{n,n-1}' \mE_{n,n-1} $.
\item $(1+\ga_n)^* =
 \Theta_{n,1}' \mP_{n,1} \Theta_{n,2}'
\mP_{n,2}\cdots \Theta_{n,n-1}' \mP_{n,n-1} $. Moreover,  for
$s=1, \ldots , n-1$, $(1+\ga_{n,s})^* =
 \Theta_{n,s}' \mP_{n,s} \Theta_{n,s+1}'
\mP_{n,s+1}\cdots \Theta_{n,n-1}' \mP_{n,n-1} $.
\end{enumerate}
\end{theorem}


\section{The groups ${\rm K}_1(\mS_n )$ and $\GL_\infty (\mS_n)$
and their generators}\label{EEE}

 In this section,  explicit generators are found for the group
 $\GL_\infty (\mS_n)$ (Theorem \ref{16Apr10}, Theorem
 \ref{17Apr10}.(1)) and it is proved that ${\rm K}_1(\mS_n) \simeq
 K^*$ (Theorem \ref{17Apr10}.(2)) modulo Theorem \ref{8May10}
 which is proved in Section \ref{THM8May10}.

The subgroup $(1+\gp_n)^*$ of the group $\mS_n^*$ is canonically
isomorphic to the group $\GL_\infty (\mS_{n-1})$ via the
isomorphism $1+\sum a_{ij}E_{ij}(n)\mapsto 1+\sum a_{ij}E_{ij}$
where $a_{ij}\in \mS_{n-1}= \bigotimes_{i=1}^{n-1}\mS_1(i)$. It is
convenient  to identify the groups $(1+\gp_n)^*$  and  $\GL_\infty
(\mS_{n-1})$ and to identify the matrix units $E_{ij}(n)$ and
$E_{ij}$, i.e. $(1+\gp_n)^*=\GL_\infty (\mS_{n-1})$ and
$E_{ij}(n)=E_{ij}$. The group $(1+\gp_n)^*$ contains the
descending chain of normal subgroups
$$(1+\gp_n)^*=(1+\gp_n)^*_1\supset \cdots \supset (1+\gp_n)^*_s\supset \cdots \supset
(1+\gp_n)^*_n=(1+F_n)^*\supset (1+\gp_n)^*_{n+1}=\{ 1\}$$ where
$(1+\gp_n)^*_s:=(1+\gp_n)^*\cap(1+\ga_{n,s})^*$. The following
lemma describes the normal subgroups $(1+\gp_n)^*_s$.

\begin{lemma}\label{a23Apr10}
$$(1+\gp_n)_s^* =\begin{cases}
(1+\sum_{|I|=s, n\in I}\gp_I)^*  & \text{if }s=1, \ldots , n-1,\\
(1+F_n)^*& \text{if }s=n.\\
\end{cases}$$
\end{lemma}

{\it Proof}. As the case $s=n$ is obvious,  we assume that $s\neq
n$. The ideal $\ga_{n,s}=\sum_{|I|=s}\gp_I$ of the algebra $\mS_n$
is the sum of idempotent ideals $\gp_I$. Therefore,
$\ga_{n,s}^2=\ga_{n,s}$. By Corollary 7.4.(3), \cite{shrekalg},
$\ga \cap \gb = \ga\gb$ for all idempotent ideals $\ga$ and $\gb$
of the algebra $\mS_n$. Since the ideals $\gp_n$ and $\ga_{n,s}$
of the algebra $\mS_n$ are idempotent, 
\begin{equation}\label{pnans}
\gp_n\cap \ga_{n,s} = \gp_n\ga_{n,s} =\sum_{|I|=s} \gp_n \gp_I=
\sum_{|I|=s, n\in I}\gp_I.
\end{equation}
Thus  $(1+\gp_n)_s^* =(1+\gp_n)^* \cap (1+\ga_{n,s})^* =
(1+\gp_n\cap \ga_{n,s})^* = (1+\sum_{|I|=s, n\in I}\gp_I)^*$.
 $\Box $

$\noindent $

For each number $s=1, \ldots , n-1$, consider the following subset
of $\mE_{n,s}$, $$\tmE_{n,s}=\prod_{|I|=s, n\in I} U_I(K)\ltimes
E_\infty (\mS_{CI}),$$ where the groups $U_I=U_I(K)$ and $E_\infty
(\mS_{CI})$ are defined in (\ref{apns3}). This is the product of
the subgroups $U_I(K)\ltimes E_\infty (\mS_{CI})$ of
$(1+\gp_n)_s^*$ in arbitrary order but which is assumed to be fixed.
Notice that $\tmE_{n,1}=U_n(K)\ltimes E_\infty (\mS_{n-1})$ where
$U_n(K)=\{ \mu_n(\l ) = \l e_n+1-e_n=\Bigl(
\begin{matrix} \l &0 \\ 0&1
\end{matrix}\Bigr)\, | \, \l \in K^*\}$ and $E_\infty (\mS_{n-1})$ is
the subgroup of $\GL_\infty (\mS_{n-1})$ generated by all the
elementary matrices.

Consider the element  $\mu_I(\l ) = \l e_I+1-e_I\in U_I$ where
$|I|=s$ and $n\in I$. Then 
\begin{equation}\label{muIG} \mu_I(\l ) = e_n(1+(\l
-1)e_{I\backslash n})+1-e_n =\Bigl(
\begin{matrix}  1+(\l -1)e_{I\backslash n} &0 \\ 0&1
\end{matrix}\Bigr)\in \GL_\infty (\mS_{n-1}).
\end{equation}

\begin{lemma}\label{a16Apr10}
$E_\infty (\mS_{n-1}) \supseteq \tmE_{n,s}$ for all $s=2, \ldots ,
n-1$.
\end{lemma}

{\it Proof}.  It is sufficient to show that the group $E_\infty
(\mS_{n-1})$ of elementary matrices  contains the groups $E_\infty
(\mS_{CI})$ and $U_I(K)$ where $|I|=s$ and $n\in I$. The group
$E_\infty (\mS_{CI})$ is generated by the elementary matrices
$u=1+aE_{\alpha\beta}(I)$ where $a\in \mS_{CI}$, $\alpha =
(\alpha_i)_{i\in I}$, $\beta = (\beta_i)_{i\in I}\in \N^I$ and
$\alpha \neq \beta$. If $\alpha_n\neq \beta_n$ then
$u=1+(a\prod_{i\in I, i\neq
n}E_{\alpha_i\beta_i}(i))E_{\alpha_n\beta_n}(n)\in E_\infty
(\mS_{n-1})$. If $\alpha_n=\beta_n$ then choose an element $\g \in
\N^I$ such that $\g_n\neq \alpha_n$, and so $\g \neq \alpha$ and
$\g \neq \beta$. Since the elements $1+E_{\alpha\g}$ and
$1+aE_{\g\beta}$ belong to the group $E_\infty (\mS_{n-1})$ (by
the previous case), so does their group commutator
$$E_\infty (\mS_{n-1})\ni [1+E_{\alpha\g}, 1+aE_{\g\beta}]=1+aE_{\alpha\beta}=u.$$
Therefore, $E_\infty (\mS_{CI})\subseteq E_\infty (\mS_{n-1})$.

It remains to show that $U_I(K)\subseteq E_\infty (\mS_{n-1})$,
i.e. $\mu_I(\l ) = 1+\l E_{00}(I)\in E_\infty (\mS_{n-1})$ for all
scalars $\l \in K\backslash \{ -1\}$. Notice that $n\in I$ and
$|I|=s\geq 2$. Choose an element, say $m\in I$, distinct from $n$.
 In the  subgroup $\GL_\infty (\mS_1(m))$ of $\GL_\infty
 (\mS_{n-1})$,  we have for all scalars $\l \in
K\backslash \{ -1\}$ the equality: 
\begin{equation}\label{AcomEi1}
\Bigl(
\begin{matrix} 1  & 0\\ -\frac{y_m}{1+\l } & 1
\end{matrix}\Bigr)\,
\Bigl(
\begin{matrix} 1 & \l x_m \\  0& 1
\end{matrix}\Bigr)\,
\Bigl(
\begin{matrix} 1 &0 \\y_m  & 1
\end{matrix}\Bigr)\,
\Bigl(
\begin{matrix} 1 & -\frac{\l}{1+\l} x_m \\  0& 1
\end{matrix}\Bigr)\,
= \Bigl(
\begin{matrix} 1+\l & 0\\ 0 & \frac{1}{1+\l }
\end{matrix}\Bigr)\,
\Bigl(
\begin{matrix} 1-\frac{\l E_{00}(m)}{1+\l} & 0\\ 0 &1
\end{matrix}\Bigr).
\end{equation}
This can be checked by direct multiplication using the equalities
$y_mx_m=1$ and  $x_my_m= 1-E_{00}(m)$  that hold in the algebra
$\mS_1(m)$. The first five matrices in the equality belong to the
group $E_\infty (\mS_1(m))$. Therefore, the last matrix $c=\Bigl(
\begin{matrix} 1-\frac{\l E_{00}(m)}{1+\l} & 0\\ 0 &1
\end{matrix}\Bigr)$ belongs to the group $E_\infty (\mS_1(m))$.
The idempotent $e:=
\begin{cases}
\prod_{i\in I\backslash \{ n,m\} }E_{00}(i)& \text{if } |I|>2,\\
1& \text{if } |I|=2, \\
\end{cases}$ determines the group monomorphism
\begin{equation}\label{taueGL}
\tau_e: \GL_\infty (\mS_1(m))=(1+\sum_{i,j\in
\N}\mS_1(m)E_{ij}(m))^*\ra \GL_\infty (\mS_{n-1})= (1+\gp_n)^*,
\;\; u\mapsto eu+1-e,
\end{equation}
that maps the group $E_\infty (\mS_1(m))$ into the group $E_\infty
(\mS_{n-1})$. Therefore,
\begin{eqnarray*}
 \tau_e(c)&=&e\bigl( E_{00}(n)(1-\frac{\l
 }{1+\l}E_{00}(m))+1-eE_{00}(n)\bigr) +1-e\\
 &=& 1-\frac{\l }{1+\l}E_{00}(I)=\mu_I(-\frac{\l }{1+\l})\in
E_\infty(\mS_{n-1})\cap U_K(I).
\end{eqnarray*}
  Since the map $\v :K\backslash
\{ -1\}\ra K\backslash \{ -1\}$, $ \l \mapsto -\frac{\l }{1+\l}$,
is a bijection ($\v^{-1} = \v$), all the elements $\mu_I(\l )$
belong to the group $E_\infty (\mS_{n-1})$. The proof of the lemma
 is complete.  $\Box $

$\noindent $

By (\ref{apns1}), there is  the group monomorphism
$$ \v_{n,s}:\frac{(1+\gp_n)^*_s}{(1+\gp_n)^*_{s+1}}\ra
\frac{(1+\ga_{n,s})^*}{(1+\ga_{n,s+1})^*}\ra
\prod_{|I|=s}(1+\bgp_I)^*= \prod_{|I|=s, n\in I}(1+\bgp_I)^*\times
\prod_{|I'|=s, n\not\in I'}(1+\bgp_{I'})^*$$ which is the
composition of two group monomorphisms. By Lemma \ref{a23Apr10},
\begin{equation}\label{ipns1}
\im (\v_{n,s})\subseteq \prod_{|I|=s, n\in I}(1+\bgp_I)^*.
\end{equation}
Recall that $(1+\bgp_I)^* = (\mX_{CI}\times U_I)\ltimes E_\infty
(L_{CI})$. Since
$\v_{n,s}(\tmE_{n,s}(1+\gp_n)^*_{s+1})=\prod_{|I|=s, n\in I}
U_I\ltimes E_\infty (L_{CI})$, we see that
$$ \v_{n,s}^{-1} (\bG_{n,s}) = \v_{n,s}^{-1} (\im (\v_{n,s})\cap\bG_{n,s}) =
\v_{n,s}^{-1} (\prod_{|I|=s, n\in I} U_I\ltimes E_\infty
(L_{CI}))=\tmE_{n,s}(1+\gp_n)^*_{s+1},
$$
and so there is the group monomorphism
$$\overline{\v}_{n,s}: (1+\gp_n)^*_s/ \tmE_{n,s}(1+\gp_n)^*_{s+1}\ra
 (1+\ga_{n,s})^* / \G_{n,s}\simeq \mX_{n,s}\simeq \mX_{n,s}'=
\prod_{|I|=s, n\in I}\mX_{CI}'\times \prod_{|I'|=s, n\not\in
I'}\mX_{CI'}'.$$ Notice that the group
$\tmE_{n,s}(1+\gp_n)^*_{s+1}$ is a normal subgroup of
$(1+\gp_n)^*_s$. For each number $s=2, \ldots , n-1$,  in the set
$\Theta_{n,s}'$ consider the exact product of cyclic groups (the
order is arbitrary) 
\begin{equation}\label{ttns}
\tTh_{n,s}:=\prod_{|J|=s+1, n\in J}\prod_{j\in J\backslash \{ n,
m(J)\}}\langle \theta_{m(J), j}(J)\rangle
\end{equation}
where $m(J)$ is the largest  element of the set $J\backslash n$.
Instead of the element $m (J)$ we can choose an arbitrary element
of the set $J\backslash n$.
 By (\ref{ipns1}),  $\im
(\overline{\v}_{n,s})\subseteq \prod_{|I|=s, n\in I}\mX_{CI}'$.
 Recall that $\im (\psi_{n,s}')= \psi_{n,s}'(\Theta_{n,s}) = \mK_{n,s}=\bigcap_{|J|=s+1}\ker (\chi_J')$,
  by Theorem \ref{24May9} and (\ref{ptijJ1}). The following
  argument is the key moment in the proof of Theorem
  \ref{16Apr10},
\begin{eqnarray*}
 \im (\overline{\v}_{n,s})&\subseteq & \im (\psi_{n,s}')\bigcap \prod_{|I|=s, n\in
 I}\mX_{CI}'=\bigcap_{|J|=s+1}\ker (\chi_J')\bigcap \prod_{|I|=s, n\in
 I}\mX_{CI}'\\
 &=&\begin{cases}
0& \text{if } s=1,\\
\prod_{|J|=s+1, n\in J} \prod_{j\in J\backslash \{ n,m(J)\}} \Z
(-(m(J), J\backslash m(J))+(j, J\backslash j))
& \text{if } s=2, \ldots , n-1,\\
\end{cases}\\
&\stackrel{{\rm by}\; (\ref{ptijJ})}{=}&\begin{cases}
0& \text{if } s=1,\\
\psi_{n,s}'(\tTh_{n,s})
& \text{if } s=2, \ldots , n-1,\\
\end{cases}\\
&=&\begin{cases}
0& \text{if } s=1,\\
\overline{\v}_{n,s}'(\tTh_{n,s}\tmE_{n,s}(1+\gp_n)^*_{s+1})
& \text{if } s=2, \ldots , n-1.\\
\end{cases}
\end{eqnarray*}
It follows that 
\begin{equation}\label{pns1}
(1+\gp_n)^*_s=\begin{cases}
\tmE_{n,1}(1+\gp_n)^*_2& \text{if } s=1,\\
\tTh_{n,s}\times_{ex}\tmE_{n,s}(1+\gp_n)^*_{s+1}
& \text{if } s=2, \ldots , n-1.\\
\end{cases}
\end{equation}

\begin{theorem}\label{16Apr10}
The group $\GL_\infty (\mS_{n-1}) = (1+\gp_n)^* $ is equal to  $
\tmE_{n,1}\tTh_{n,2}\tmE_{n,2}\cdots \tTh_{n,n-1}\tmE_{n,n-1}$.
 Moreover,
$$ (1+\gp_n)_s^*=\begin{cases}
\tmE_{n,1}\tTh_{n,2}\tmE_{n,2}\cdots \tTh_{n,n-1}\tmE_{n,n-1}& \text{if }s=1,\\
\tTh_{n,s}\tmE_{n,s}\cdots \tTh_{n,n-1}\tmE_{n,n-1}& \text{if }s=2, \ldots, n-1,\\
(1+F_n)^* &  \text{if }s=n.\\
\end{cases}$$
\end{theorem}

{\it Proof}. By Proposition 3.10, \cite{Snaut}, we have the
inclusion  $(1+\gp_n)_n^*=(1+F_n)^* \subseteq \tmE_{n,n-1}$. Now,
the theorem follows from (\ref{pns1}). $\Box $

$\noindent $

For each subset $J$ of the set $\{ 1, \ldots , n\}$ such that
$n\in J$ and $|J|\geq 3$, and for each pair of distinct elements
$i$ and $j$ of the set  $J\backslash n$, the unit $\th_{ij}(J)\in
\mS_n^*$ can be written as follows
\begin{eqnarray*}
  \th_{ij}(J)&=&(y_ie_{J\backslash i}e_n+1-e_n+e_n(1-e_{J\backslash
i}))(x_je_{J\backslash j}e_n+1-e_n+e_n(1-e_{J\backslash j}))\\
 &=&e_n(y_ie_{J\backslash i}+1-e_{J\backslash
i})(x_je_{J\backslash j}+1-e_{J\backslash j})+1-e_n\\
&=&e_n\th_{ij}(J\backslash n)+1-e_n
\end{eqnarray*}
where $e_n:= E_{00}(n)$, $e_{J\backslash i}:=\prod_{k\in
J\backslash i} E_{00}(k)$ and $e_{J\backslash j}:=\prod_{k\in
J\backslash j} E_{00}(k)$. Therefore, the unit $\th_{ij}(J)$, as
an element of the group $\GL_\infty (\mS_{n-1})$,  is the matrix
\begin{equation}\label{tijJm}
\th_{ij}(J)=\Bigl(
\begin{matrix}  \th_{ij}(J\backslash n) &0 \\ 0&1
\end{matrix}\Bigr)\in \GL_\infty (\mS_{n-1})
\end{equation}
where $\th_{ij}(J\backslash n)\in \mS_{n-1}^*$.

{\bf The determinant $\bdet$ on $\GL_\infty (\mS_{n-1})$}. The
algebra epimorphism $\mS_{n-1}\ra \mS_{n-1}/\ga_{n-1}=L_{n-1}$,
$a\mapsto \overline{a} := a+\ga_{n-1}$, yields the group
homomorphisms $\GL_\infty (\mS_{n-1})\ra \GL_\infty (L_{n-1})$,
$u\mapsto \bu $, and $\bdet : \GL_\infty (\mS_{n-1})\ra
\GL_{\infty} (L_{n-1})\stackrel{\det}{\ra} L_{n-1}^*$. Clearly,
$\bdet (E_\infty (\mS_{n-1}))=1$, $\bdet (\tTh_{n,s})=1$ for all
$s=2, \ldots , n-1$, and $\bdet (U_n(K))=K^*$ since $\bdet
(\mu_n(\l )) = \l$ for all $\l \in K^*$. By Theorem \ref{16Apr10}
and Lemma \ref{a16Apr10}, $\GL_\infty (\mS_{n-1}) = U_n(K)
\tTh_{n,2}\cdots \tTh_{n,n-1}E_\infty (\mS_{n-1})$,  since $E_\infty
(\mS_{n-1})$ is a normal subgroup of $\GL_\infty (\mS_{n-1})$. It
follows that the image of the map $\bdet$ is $K^*$, i.e. we have
the group epimorphism 
\begin{equation}\label{bn1det}
\bdet :\GL_\infty (\mS_{n-1})\ra K^*, \;\; u\mapsto \det
(\overline{u}),
\end{equation}
and 
\begin{equation}\label{cn1det}
\GL_\infty (\mS_{n-1})= U_n(K)\ltimes \ker (\bdet ), \;\;
\SL_\infty (\mS_{n-1}):=\ker (\bdet ) =\tTh_{n,2}\cdots
\tTh_{n,n-1}E_\infty (\mS_{n-1}).
\end{equation}
\begin{theorem}\label{8May10}
$\tTh_{n,s}\subseteq E_\infty (\mS_{n-1})$ for all $s=2, \ldots ,
n-1$.\end{theorem}


The proof of Theorem \ref{8May10} is not easy and is given in
Section \ref{THM8May10}.

\begin{theorem}\label{17Apr10}
\begin{enumerate}
\item $\GL_\infty (\mS_{n-1})=U_n(K)\ltimes E_\infty (\mS_{n-1})$
 and $\SL_\infty (\mS_{n-1}) = E_\infty (\mS_{n-1})$ where $U_n(K) =\{ \mu_n(\l ) := 1+(\l -1)E_{00}(n)\, |
\, \l\in K^*\}$. So, each element $a\in \GL_\infty (\mS_{n-1})$ is
the unique product $a= \mu_n(\l ) e$ where $\l = \bdet (a)$ and
$e:= \mu_n(\bdet (a))^{-1} a\in E_\infty (\mS_{n-1})$. \item ${\rm
K}_1(\mS_n)\simeq K^*$ for all $n\geq 1$.
\end{enumerate}
\end{theorem}

{\it Proof}. The theorem follows from Theorem \ref{8May10} and
(\ref{cn1det}). $\Box $

$\noindent $

The number of generators $\th_{\max (J), j}(J)$ in the block
$\tTh_{n+1,2}\cdots \tTh_{n+1,n}$  for the group $\GL_\infty
(\mS_n)= U_{n+1}(K)\ltimes \tTh_{n+1,2}\cdots \tTh_{n+1,n}E_\infty
(\mS_n)$ is $\sum_{s=2}^n {n\choose s} (s-1) = (n-2)2^{n-1}+1$ as
the next lemma shows.

\begin{lemma}\label{a17Apr10}
For each natural number $n\geq 2$, $\sum_{s=2}^n {n\choose s}
(s-1) = (n-2)2^{n-1}+1$.
\end{lemma}

{\it Proof}. Taking the derivative of the polynomial $(1+x)^n =
\sum_{s=0}^n {n\choose s} x^s$, we have the equality $n(1+x)^{n-1}
= \sum_{s=1}^n{n\choose s} sx^{s-1}$. Then taking the difference
of both equalities at $x=1$, we obtain the result: $\sum_{s=2}^n
{n\choose s}(s-1)-1= n2^{n-1}-2^n = (n-2)2^{n-1}$. $\Box $



\section{Proof of Theorem
\ref{8May10}}\label{THM8May10}

The whole   section is a proof of Theorem \ref{8May10}. The proof
is constructive (but slightly technical) and split into a series
of lemmas that produce more and more sophisticated elementary
matrices in $E_\infty (\mS_{n-1})$. These elementary matrices are
used to show that the elements of the sets $\tTh_{n,s}$ are
elementary matrices (Propositions \ref{7May10} and \ref{9May10}).

\begin{lemma}\label{a8May10}
Let $D$ be a division ring and let $\L =D\oplus De$ be a ring over $D$
such that $e^2=e$ and $de=ed$ for all $d\in D$. Then
\begin{enumerate}
\item the group of units $\L^*$ of the ring $\L$ is the
semi-direct product $D^*\ltimes \G$ of the group of units $D^*$ of
 the ring $D$ and the subgroup $\G := \{ 1+\l e\, | \, \l \in
 D\backslash \{  -1\}\}$ of $\L^*$.
\item $(1+\l e)^{-1}=1-\frac{\l}{1+\l}e$ for all elements $\l \in
 D\backslash \{  -1\}$.
\item The map $\phi :D\backslash \{  -1\}\ra D\backslash \{ -1\}$,
 $\l\mapsto -\frac{\l}{1+\l}$, is a bijection with $\phi^{-1} =
\phi$.
\item $(1-2e)^{-1}=1-2e$.
\end{enumerate}
\end{lemma}

{\it Proof}.  Straightforward. $\Box $

We are interested in the rings $\L$ and their groups of units,
since the algebra $K+M_\infty (\mS_{n-1})$ of infinite dimensional
matrices over the algebra $\mS_{n-1}$ contains plenty of them and
as a  result the group $\GL_\infty (\mS_{n-1})$ contains their
groups of units.

\begin{lemma}\label{b8May10}
Let $\mS_1(\L )=\L \langle x,y\, | \, yx=1\rangle$ be the algebra
$\mS_1$ over the ring $\L$ from Lemma \ref{a8May10}. Then,  for
each element $\l \in D\backslash \{ -1\}$,
\begin{equation}\label{BcomE}
\Bigl(
\begin{matrix} 1  & 0\\ -\frac{y}{1+\l e} & 1
\end{matrix}\Bigr)\,
\Bigl(
\begin{matrix} 1 & \l ex \\  0& 1
\end{matrix}\Bigr)\,
\Bigl(
\begin{matrix} 1 &0 \\y  & 1
\end{matrix}\Bigr)\,
\Bigl(
\begin{matrix} 1 & -\frac{\l e}{1+\l e} x \\  0& 1
\end{matrix}\Bigr)\,
= \Bigl(
\begin{matrix} 1+\l e& 0\\ 0 & \frac{1}{1+\l e}
\end{matrix}\Bigr)\,
\Bigl(
\begin{matrix} 1-\frac{\l }{1+\l}eE_{00} & 0\\ 0 &1
\end{matrix}\Bigr).
\end{equation}
where $E_{00}:=1-xy$ (the element $1+\l e$ is a unit of the
algebra $\mS_1(\L )$, by Lemma \ref{a8May10}).
\end{lemma}

{\it Proof}. The RHS of the equality (\ref{BcomE}) is the product
of four matrices, say $A_1\cdots A_4$.
$$ A_1A_2A_3=\Bigl(
\begin{matrix} 1  & \l ex\\ -\frac{y}{1+\l e} & 1-\frac{\l e}{1+\l e}
\end{matrix}\Bigr)\,
\Bigl(
\begin{matrix} 1 &0 \\y  & 1
\end{matrix}\Bigr)=\Bigl(
\begin{matrix} 1+\l exy  & \l ex\\ 0 & \frac{1}{1+\l e}
\end{matrix}\Bigr), \, \;\;\; A_1\cdots A_4=\Bigl(
\begin{matrix} 1+\l exy  & 0\\ 0 & \frac{1}{1+\l e}
\end{matrix}\Bigr), $$
since $(1+\l exy) (-\frac{\l e}{1+\l e} x) +\l ex = -\frac{\l
e}{1+\l e}(1+\l e) x+\l ex= 0$. Now, 
\begin{equation}\label{ABcomE}
A_1\cdots A_4= \Bigl(
\begin{matrix} 1+\l e & 0\\ 0 & \frac{1}{1+\l e}
\end{matrix}\Bigr)\,
\Bigl(
\begin{matrix} 1-\frac{\l e}{1+\l e}E_{00}  & 0\\ 0 & 1
\end{matrix}\Bigr)
\end{equation}
since $(1+\l e) (1-\frac{\l e}{1+\l e}E_{00}) =1+ \l e
(1-E_{00})=1+\l exy$. Finally, the equality (\ref{BcomE}) follows
from Lemma \ref{a8May10}.(2), $\frac{\l e}{1+\l e}=\l e(1-\frac{\l
}{1+\l }e)=\l (1-\frac{\l }{1+\l })e=\frac{\l }{1+\l }e$. $\Box $

$\noindent $

For each ring $R$ and a natural number $m\geq 1$, $E_n(R)$ is the
subgroup of $\GL_n(R)$ generated by all  elementary matrices.

\begin{lemma}\label{a7May10}
\begin{enumerate}
\item $\Bigl(
\begin{matrix}  y &0 \\ E_{00} &x
\end{matrix}\Bigr)\in E_2(\mS_1)$ where $E_{00}:= 1-xy$.
\item $\Bigl(
\begin{matrix} x  & E_{00}\\0 &y
\end{matrix}\Bigr)=\Bigl(
\begin{matrix}  y &0 \\ E_{00}&x
\end{matrix}\Bigr)^{-1}\in E_2(\mS_1)$.
\end{enumerate}
\end{lemma}

{\it Proof}. 1. Using the equalities $yx=1$ and $E_{00}x=0$,  we can
easily check that 
\begin{equation}\label{1AcomE}
\Bigl(
\begin{matrix} 1  & 0\\ 1 & 1
\end{matrix}\Bigr)\,
\Bigl(
\begin{matrix} 1 & -1\\  0& 1
\end{matrix}\Bigr)\,
\Bigl(
\begin{matrix} 1 &0 \\1-x  & 1
\end{matrix}\Bigr)\,
\Bigl(
\begin{matrix} y & 0\\  E_{00}& x
\end{matrix}\Bigr)\,\Bigl(
\begin{matrix} 1 & x\\ 0 & 1
\end{matrix}\Bigr)\,\Bigl(
\begin{matrix} 1 & 0\\ -y & 1
\end{matrix}\Bigr)\,
= \Bigl(
\begin{matrix} 1- 2E_{00} & 0\\ 0 &1
\end{matrix}\Bigr).
\end{equation}
By (\ref{AcomEi1}), the RHS is an element of the  group
$E_2(\mS_1)$ since $\frac{-2}{1+(-2)}=2$, and so statement 1
holds.

2. It is obvious. $\Box $


Let $R$ be a ring and $u$ be its unit. The $2\times2$ matrix
$\Bigl(
\begin{matrix} y & 0\\ uE_{00} &x
\end{matrix}\Bigr)\in M_2(\mS_1(R))$ is invertible where $E_{00}:=1-xy$. Moreover,
\begin{equation}\label{yuE}
\Bigl(
\begin{matrix} y & 0\\ uE_{00} &x
\end{matrix}\Bigr)^{-1}= \Bigl(
\begin{matrix} x & u^{-1}E_{00}\\ 0 &y
\end{matrix}\Bigr).
\end{equation}
\begin{lemma}\label{c8May10}
Let the ring $\L$ be as in Lemma \ref{a8May10}. Then, for each
element $\l \in D\backslash \{ -1\}$, $$\Bigl(
\begin{matrix} y  &0 \\(1+\l e)E_{00} &x
\end{matrix}\Bigr)\in E_2(\mS_1(\L )) \;\; {\rm and }\;\; \Bigl(
\begin{matrix} x  &(1+\l e)^{-1}E_{00} \\0 &y
\end{matrix}\Bigr)= \Bigl(
\begin{matrix} y  &0 \\ (1+\l e)E_{00}& x
\end{matrix}\Bigr)^{-1}\in E_2(\mS_1(\L ))$$ where $(1+\l e)^{-1} =
1-\frac{\l}{1+\l}e$ (by Lemma \ref{a8May10}.(2)).

\end{lemma}

{\it Proof}. It suffices to prove the first inclusion since
then the equality and the second inclusion follow from
(\ref{yuE}). Using the equalities $yx=1$ and $E_{00}x=0$ we can
check that 
\begin{equation}\label{1yuE}
\Bigl(
\begin{matrix} 1  & 0\\ 1 & 1
\end{matrix}\Bigr)\,
\Bigl(
\begin{matrix} 1 & -1\\  0& 1
\end{matrix}\Bigr)\,
\Bigl(
\begin{matrix} 1 &0 \\1-x  & 1
\end{matrix}\Bigr)\,
\Bigl(
\begin{matrix} y & 0\\ (1+\l e) E_{00}& x
\end{matrix}\Bigr)\,\Bigl(
\begin{matrix} 1 & x\\ 0 & 1
\end{matrix}\Bigr)\,\Bigl(
\begin{matrix} 1 & 0\\ -y & 1
\end{matrix}\Bigr)\,
= \Bigl(
\begin{matrix} 1- (2+\l e)E_{00} & 0\\ 0 &1
\end{matrix}\Bigr),
\end{equation}
$$\Bigl(
\begin{matrix} 1- (2+\l e)E_{00} & 0\\ 0 &1
\end{matrix}\Bigr)=\Bigl(
\begin{matrix} 1- 2E_{00} & 0\\ 0 &1
\end{matrix}\Bigr)\, \Bigl(
\begin{matrix} 1+\l eE_{00} & 0\\ 0 &1
\end{matrix}\Bigr).$$
By (\ref{1AcomE}), $\Bigl(
\begin{smallmatrix} 1- 2E_{00} & 0\\ 0 &1
\end{smallmatrix}\Bigr)\in E_2(\mS_1)$, and then  by (\ref{BcomE}), $ \Bigl(
\begin{smallmatrix} 1+\l eE_{00} & 0\\ 0 &1
\end{smallmatrix}\Bigr)\in E_2(\mS_1(\L ))$ since $\l \in D\backslash \{ -1\}$.
Therefore, $\Bigl(
\begin{smallmatrix} y & 0\\ (1+\l e) E_{00}& x
\end{smallmatrix}\Bigr)\in E_2(\mS_1(\L ))$, by (\ref{1yuE}). $\Box $


\begin{lemma}\label{b7May10}
$\Bigl(
\begin{matrix} 1+(y_2-1)x_1y_1  &0 \\ e_2y_1& x_2
\end{matrix}\Bigr)\in E_2(\mS_2)$ where $e_2:=
E_{00}(2)=1-x_2y_2$.
\end{lemma}

{\it Proof}. The statement follows from the equality
\begin{equation}\label{3AcomE}
\Bigl(
\begin{matrix} 1  & 0\\ -x_2y_1 & 1
\end{matrix}\Bigr)\,
\Bigl(
\begin{matrix} 1 & (y_2-1)x_1\\  0& 1
\end{matrix}\Bigr)\,
\Bigl(
\begin{matrix} 1 &0 \\y_1  & 1
\end{matrix}\Bigr)\,
\Bigl(
\begin{matrix} 1 & -(y_2-1)x_1\\ 0& 1
\end{matrix}\Bigr)\,\Bigl(
\begin{matrix} 1 & (y_2-1)(1-x_2)x_1\\ 0 & 1
\end{matrix}\Bigr)
= \Bigl(
\begin{matrix} 1+ (y_2-1)x_1y_1 & 0\\ e_2y_1 &x_2
\end{matrix}\Bigr)
\end{equation}
which can be checked directly using the equalities $y_ix_i=1$,
$x_iy_i=1-e_i$, $y_ie_i=0$  and $e_ix_i=0$ where $e_i:=
E_{00}(i)$. The RHS of the equality (\ref{3AcomE}) is the product
of five matrices $A_1\cdots A_5$.
$$ A_1A_2A_3=\Bigl(
\begin{matrix}  1 &(y_2-1)x_1 \\ -x_2y_1& 1-x_2(y_2-1)
\end{matrix}\Bigr)\, \Bigl(
\begin{matrix} 1  &0 \\ y_1&1
\end{matrix}\Bigr)=\Bigl(
\begin{matrix}  1+(y_2-1)x_1y_1 &(y_2-1)x_1 \\ e_2y_1& 1-x_2(y_2-1)
\end{matrix}\Bigr)$$
since $-x_2y_1+(1-x_2y_2)y_1+x_2y_1= e_2y_1$. Now,
$$ A_1\cdots A_4=\Bigl(
\begin{matrix}  1+(y_2-1)x_1y_1 &-(y_2-1)^2x_1 \\ e_2y_1& 1-(x_2+e_2)(y_2-1)
\end{matrix}\Bigr)$$
since $-y_1x_1e_2(y_2-1) +1-x_2(y_2-1) = 1-(x_2+e_2)(y_2-1)$.
Finally, $A_1\cdots A_5=\Bigl(
\begin{matrix} 1+(y_2-1)x_1y_1  &a \\ e_2y_1& b
\end{matrix}\Bigr)$ where
\begin{eqnarray*}
 a&=&(1+(y_2-1)x_1y_1)(y_2-1) (1-x_2)x_1-(y_2-1)^2x_1 \\
 &=& (x_1+(y_2-1)x_1)(y_2-1)(1-x_2)-(y_2-1)^2x_1\\
 &=&x_1(y_2-1) (y_2-1)-(y_2-1)^2x_1=0,\\
 b&=& 1-(x_2+e_2)(y_2-1) +e_2y_1(y_2-1)(1-x_2) x_1 \\
 &=& 1-x_2(y_2-1) -e_2(y_2-1) +e_2(y_2-1) - e_2(1-x_2)\\
 &=& 1-x_2y_2+x_2-e_2= x_2. \;\; \Box
\end{eqnarray*}

\begin{proposition}\label{7May10}
$\Bigl(
\begin{matrix} \th_{12}  & 0\\ 0&1
\end{matrix}\Bigr)\in E_2(\mS_2)$ where $\th_{12} = \th_{12}(\{ 1,2\}) =
(1+(y_1-1)e_2)(1+(x_2-1)e_1)$, $e_1= E_{00}(1)$ and $e_2=
E_{00}(2)$.
\end{proposition}

{\it Proof}. By Lemma \ref{a7May10}, $\Bigl(
\begin{matrix}  x_2 &e_2 \\ 0&y_2
\end{matrix}\Bigr)\in E_2(\mS_1(2))\subseteq E_2(\mS_2)$. Then, by
Lemma \ref{b7May10}, 
\begin{equation}\label{2AcomE}
 E_2(\mS_2)\ni \Bigl(
\begin{matrix}  x_2 &e_2 \\ 0&y_2
\end{matrix}\Bigr)\, \Bigl(
\begin{matrix} 1+(y_2-1)x_1y_1  &0 \\ e_2y_1& x_2
\end{matrix}\Bigr)=\Bigl(
\begin{matrix} \th_{12}  & 0\\ 0&1
\end{matrix}\Bigr).
\end{equation}
Indeed, let $a$ be the $(1,1)$-entry of the product, then
\begin{eqnarray*}
 a&=& x_2(1+(y_2-1)x_1y_1) +e_2^2y_1=
 x_2+(x_2y_2-x_2)x_1y_1+e_2y_1 \\
 &=&x_2e_1+(1-e_2) (1-e_1) +e_2y_1= 1+(x_2-1)e_1+(y_1-1) e_2+e_1e_2 \\
 &=& (1+(y_1-1)e_2)(1+(x_2-1)e_1) = \th_{12}
\end{eqnarray*}
since $(y_1-1)e_2\cdot (x_2-1) e_1= (y_1-1) e_1\cdot e_2(x_2-1) =
(-e_1)\cdot (-e_2)= e_1e_2$. $\Box $

\begin{lemma}\label{a9May10}
Let $J=\{ 1, \ldots , m\}$ where $m\geq 3$, and let $I=J\backslash
\{ 1,2\}$. Then $$\Bigl(
\begin{matrix} 1+(y_2-1)x_1y_1e_I  &0 \\ e_2y_1e_I& 1+(x_2-1)e_I
\end{matrix}\Bigr)\in E_2(\mS_2(K\oplus Ke_I))$$ where
$e_2:=E_{00}(2)$ and $e_I:= \prod_{k\in I}E_{00}(k)$.
\end{lemma}

{\it Proof}. The statement follows from the equality
\begin{equation}\label{1BcomE}
\Bigl(
\begin{smallmatrix} 1  & 0\\ -x_2y_2e_I & 1
\end{smallmatrix}\Bigr)\,
\Bigl(
\begin{smallmatrix} 1 & (y_2-1)x_1\\  0& 1
\end{smallmatrix}\Bigr)\,
\Bigl(
\begin{smallmatrix} 1 &0 \\y_1e_I  & 1
\end{smallmatrix}\Bigr)\,
\Bigl(
\begin{smallmatrix} 1 & -(y_2-1)x_1\\ 0& 1
\end{smallmatrix}\Bigr)\,\Bigl(
\begin{smallmatrix} 1 & (y_2-1)(1-x_2)x_1e_I\\ 0 & 1
\end{smallmatrix}\Bigr)
= \Bigl(
\begin{smallmatrix} 1+ (y_2-1)x_1y_1e_I & 0\\ e_2y_1e_I &1+(x_2-1)e_I
\end{smallmatrix}\Bigr).
\end{equation}
The equality can be written shortly as $A_1\cdots A_5=A$.
$$ A_2A_3A_4=\Bigl(
\begin{matrix}  1+(y_2-1)x_1y_1e_I & (y_2-1)x_1\\ y_1e_I&1
\end{matrix}\Bigr)\, \Bigl(
\begin{matrix} 1  & -(y_2-1)x_1\\ 0&1
\end{matrix}\Bigr)=\Bigl(
\begin{matrix}  1+(y_2-1)x_1y_1e_I & -(y_2-1)^2x_1e_I\\
y_1e_I&1-(y_2-1)e_I
\end{matrix}\Bigr)
$$
where we have used the fact that $y_1x_1=1$.
$$ A_1\cdots A_4=\Bigl(
\begin{matrix}  1+(y_2-1)x_1y_1e_I & -(y_2-1)^2x_1e_I\\
e_2y_1e_I&1-(x_2+e_2)(y_2-1)e_I
\end{matrix}\Bigr).
$$
In more detail, let $(\alpha, \beta )$ be the second row of the
product. Using the fact that $y_1x_1=1$ and $e_I^2=e_I$, we see
that
\begin{eqnarray*}
 \alpha &=& -x_2y_1e_I(1+(y_2-1)x_1y_1e_I) +y_1e_I= (-x_2(y_1+(y_2-1)y_1)+y_1) e_I \\
 &=& (1-x_2y_2)y_1e_I= e_2y_1e_I,
\end{eqnarray*}
\begin{eqnarray*}
 \beta &=& x_2y_1e_I(y_2-1)^2x_1e_I+1-(y_2-1)e_I= 1+(x_2y_2-x_2-1) (y_2-1)e_I\\
 &=& 1-(x_2+e_2)(y_2-1)e_I.
\end{eqnarray*}
Finally, $A_1\cdots A_5=\Bigl(
\begin{matrix}  1+(y_2-1)x_1y_1e_I & a'\\
e_2y_1e_I&b'
\end{matrix}\Bigr)$ where (below, we use the fact  that $a=0$ and
$b=x_2$, see the proof of Lemma \ref{b7May10})
\begin{eqnarray*}
 a' &=& (1+(y_2-1)x_1y_1e_I)(y_2-1)
 (x_2-1)x_1e_I-(y_2-1)^2x_1e_I\\
 &=& \Bigl( (1+(y_2-1)x_1y_1)(y_2-1)
 (x_2-1)x_1-(y_2-1)^2x_1   \Bigr)e_I=a\cdot e_I= 0\cdot e_I=0,\\
\beta &=& 1-(x_2+e_2)(y_2-1)e_I+e_2y_1(y_2-1)(1-x_2)x_1e_I\\
 &=& 1+\Bigl(
 -1+1-(x_2+e_2)(y_2-1)+e_2y_1(y_2-1)(1-x_2)x_1\Bigr)e_I\\
 &=&1+(-1+b)e_I = 1+(x_2-1)e_I.
\end{eqnarray*}
The proof of the lemma is complete.  $\Box $

$\noindent $

Let $J=\{ 1,2, \ldots , m\}$ and $m\geq 3$. By multiplying out,
the element $\th_{12}(J)= (1+(y_1-1)e_{J\backslash
1})(1+(x_2-1)e_{J\backslash 2})\in \mS_m^*$ can be written as the
sum 
\begin{equation}\label{t12s}
\th_{12}(J)= x_2e_1e_I+(1-e_1e_I)(1-e_2e_I)+y_1e_2e_I
\end{equation}
where $I:=J\backslash \{ 1,2\}$.

\begin{proposition}\label{9May10}
Let $J=\{ 1,2, \ldots , m\}$
 and $m\geq 3$. Then
$\Bigl(
\begin{matrix} \th_{12}(J)  & 0\\ 0&1
\end{matrix}\Bigr)\in E_2(\mS_m)$.
\end{proposition}

{\it Proof}. We keep the notation of Lemma \ref{a9May10}.  By
Lemma \ref{a7May10}.(2) and Lemma \ref{a9May10}, the product of
the following two elementary matrices is also an elementary
matrix, 
\begin{equation}\label{2BcomE}
 E_2(\mS_2)\ni \Bigl(
\begin{matrix}  x_2 &e_2 \\ 0&y_2
\end{matrix}\Bigr)\, \Bigl(
\begin{matrix} 1+(y_2-1)x_1y_1e_I  &0 \\ e_2y_1e_I& 1+(x_2-1)e_I
\end{matrix}\Bigr)=\Bigl(
\begin{matrix} \th_{12}(J)+(x_2-1)(1-e_I)  & e_2(1-e_I)\\
0&e_I+(1-e_I)y_2
\end{matrix}\Bigr).
\end{equation}
Indeed, the LHS is a matrix of type $\Bigl(
\begin{matrix}  \alpha & \g \\0 &\beta
\end{matrix}\Bigr)$ (since $y_2e_2=0$) where
\begin{eqnarray*}
 \alpha &=& x_2(1+(y_2-1)x_1y_1e_I) +e_2y_1e_I = x_2+(1-e_2-x_2)(1-e_1)e_I+e_2y_1e_I\\
&=& x_2-x_2(1-e_1)e_I+
(1-e_1)(1-e_2)e_I+y_1e_2e_I\\
&=&x_2(1-e_I)+\Bigl(
x_2e_1e_I+(1-e_2e_I)(1-e_1e_I)+y_1e_2e_I\Bigr) + (1-e_1)
(1-e_2)e_I-(1-e_1e_I)(1-e_2e_I)\\
&\stackrel{{\rm by}\; (\ref{t12s})}{=}&
x_2(1-e_I)+\th_{12}(J)+e_I-e_1e_I-e_2e_I+e_J-1+e_1e_I+e_2e_I-e_J\\
&=& \th_{12}(J)+(x_2-1)(1-e_I),\\
 \beta&=& y_2(1+(x_2-1)e_I)= y_2+(1-y_2)e_I=e_I+(1-e_I)y_2,  \\
 \g &=& e_2(1+(x_2-1)e_I)=e_2(1-e_I),
\end{eqnarray*}
since $e_2x_2=0$. By  (\ref{t12s}), 
\begin{equation}\label{t12eI}
\th_{12}(J) (1-e_I)=1-e_I.
\end{equation}
Using (\ref{t12eI}), the RHS of (\ref{2BcomE}) is equal to the
product of two matrices
$$\Bigl(
\begin{matrix} \th_{12}(J)+(x_2-1)(1-e_I)  & e_2(1-e_I) \\0 &
e_I+(1-e_I)y_2
\end{matrix}\Bigr)= \Bigl(
\begin{matrix} \th_{12}(J)  &0 \\0 &1
\end{matrix}\Bigr)\, \Bigl(
\begin{matrix} 1+(x_2-1)(1-e_I)  &e_2(1-e_I) \\ 0&e_I+(1-e_I)y_2
\end{matrix}\Bigr).$$
In order to finish the proof of the proposition, it suffices to
show that the last matrix is elementary. This follows from the
next two equalities, as the last two matrices in the equality
(\ref{5BcomE}) belong to the group $E_2(\mS_m)$, by Lemma
\ref{a16Apr10}.

\begin{equation}\label{4BcomE}
\Bigl(
\begin{smallmatrix} 1  & -(x_2-1+2e_2)(1-e_I)\\0 &1
\end{smallmatrix}\Bigr)\, \Bigl(
\begin{smallmatrix}  1+(x_2-1)(1-e_I) & e_2(1-e_I)\\ 0& e_I+(1-e_I)y_2
\end{smallmatrix}\Bigr)\, \Bigl(
\begin{smallmatrix} 1  &0 \\x_2 &1
\end{smallmatrix}\Bigr)\, \Bigl(
\begin{smallmatrix}  1 &(1-y_2)(1-e_I) \\0 &1
\end{smallmatrix}\Bigr)\, \Bigl(
\begin{smallmatrix} 1  & 0\\ -1-(x_2-1)e_I& 1
\end{smallmatrix}\Bigr)=\Bigl(
\begin{smallmatrix} 1-2e_2(1-e_I)  & 0\\ 0&1
\end{smallmatrix}\Bigr),
\end{equation}
\begin{equation}\label{5BcomE}
\Bigl(
\begin{smallmatrix} 1-2e_2(1-e_I)  & 0\\ 0&1
\end{smallmatrix}\Bigr)=\Bigl(
\begin{smallmatrix} 1-2e_2  & 0\\ 0&1
\end{smallmatrix}\Bigr)\, \Bigl(
\begin{smallmatrix} 1-2e_2e_I  & 0\\ 0&1
\end{smallmatrix}\Bigr).
\end{equation}
The equality (\ref{5BcomE}) is obvious, and the equality
(\ref{4BcomE}) can be written in the form $A_1\cdots A_5=A$. Using
the identities $e_2x_2=0$, $y_2x_2=1$, $e_I^2=e_I$ and
$(1-e_I)^2=1-e_I$, we see that
$$A_2A_3A_4=\Bigl(
\begin{smallmatrix} 1+(x_2-1)(1-e_I)  &e_2(1-e_I) \\ 1+(x_2-1)e_I& e_I+(1-e_I)y_2
\end{smallmatrix}\Bigr)\, \Bigl(
\begin{smallmatrix} 1  &(1-y_2)(1-e_I) \\0 &1
\end{smallmatrix}\Bigr)=\Bigl(
\begin{smallmatrix} 1+(x_2-1)(1-e_I)  & (x_2-1+2e_2)(1-e_I)\\
1+(x_2-1)e_I&1
\end{smallmatrix}\Bigr).
$$
In more detail, let $(u,v)^t$ be the second column of the product
of the two matrices in the middle. Then
\begin{eqnarray*}
 u&=& (1+(x_2-1)(1-e_I))(1-y_2)(1-e_I)+e_2(1-e_I)=
(x_2(1-y_2)+e_2)(1-e_I)  \\
 &= & (x_2-(1-e_2)+e_2) (1-e_I)= (x_2-1+2e_2) (1-e_I),\\
 v&=&(1+(x_2-1)e_I) (1-y_2)(1-e_I)+e_I+(1-e_I)y_2\\
 &=& (1-y_2)(1-e_I)+e_I+(1-e_I)y_2=1.
\end{eqnarray*}
Finally,
$$A_2\cdots A_5=\Bigl(
\begin{matrix}  1-2e_2(1-e_I) & (x_2-1+2e_2)(1-e_I) \\0 &1
\end{matrix}\Bigr)$$
since $ 1+(x_2-1)(1-e_I)-(x_2-1+2e_2)(1-e_I)(1+(x_2-1)e_I)=
1+(x_2-1-x_2+1-2e_2)(1-e_I) = 1-2e_2(1-e_I)$. Now, (\ref{4BcomE})
is obvious. The proof of the proposition is complete. $\Box $

$\noindent $

{\bf Proof of Theorem \ref{8May10}}. Notice that $\mS_{n-1}\simeq
\mS_1^{\t (n-1)}$ and the symmetric group $S_{n-1}$ is a subgroup
of the group of automorphisms of the algebra $\mS_{n-1}$ (it acts
by permuting the tensor components). Then, the matrix $\Bigl(
\begin{matrix}  \th_{ij}(J) & 0 \\0 &1
\end{matrix}\Bigr)$  (where $J\subseteq \{ 1, \ldots , n-1\}$ with $|J|\geq 2$) is elementary by Proposition \ref{7May10}
(when $|J|=2$) and Proposition \ref{9May10} (when $|J|>2$).  Now,
Theorem \ref{8May10} is obvious. $\Box$



\section{The groups $\K1 (\mS_n, \gp)$ and
$\GL_\infty (\mS_n, \gp )$ and their
generators}\label{K1GPSN}

In this section, explicit generators are found for the group
$\GL_\infty (\mS_{n-1}, \gp)$ where $\gp$ is an arbitrary nonzero
 idempotent prime ideal of the algebra $\mS_{n-1}$ and it is proved that
 ${\rm K}_1(\mS_{n-1}, \gp )\simeq \Z^{m\choose 2}\times K^{*m}$
 (Theorem \ref{13May10}) where $m$ is the height of the ideal
 $\gp$.

For a ring $A$ and an ideal $\ga$ of $A$, the normal subgroup of
$\GL_\infty (A)$, $$\GL_\infty (A, \ga ):= \ker ( \GL_\infty (A)
\ra \GL_\infty (A/\ga )),$$ is called the {\em congruence group}
of level $\ga$. The normal subgroup $E_\infty (A, \ga )$ of $E_\infty
(A)$ which is generated by all the $\ga$-{\em elementary
matrices} ($1+a E_{ij}$, $a\in \ga$, $i\neq j$) is a normal
subgroup of $\GL_\infty (A)$. Moreover, $[ \GL_\infty (A),
\GL_\infty (A, \ga ) ] = E_\infty (A, \ga )$
\cite{Bass-book-K-theory}, and so the ${\rm K}_1$-group
$$ {\rm K}_1(A, \ga ) := \GL_\infty (A, \ga ) / E_\infty (A, \ga
)$$ is abelian. Let $E_\infty'(A, \ga )$ be the subgroup of $E_\infty (A)$ generated by all the $\ga$-elementary matrices. Then $E_\infty'(A, \ga )\subseteq E_\infty(A, \ga )\subseteq E_\infty (A)$.

We keep the notation of the previous sections. Recall that we
identified the groups $(1+\gp_n)^*$ and $\GL_\infty (\mS_{n-1})$.
Each nonzero idempotent prime ideal $\gp$ of the algebra
$\mS_{n-1}$ is a {\em unique} sum (up to order) of distinct height
one prime ideals $\gp = \gp_{i_1}+\cdots + \gp_{i_m}$ and ${\rm
ht} (\gp ) = m$ where ${\rm ht}$ stands for the {\em height} of an
ideal, Corollary 4.8, \cite{shrekalg}. The set $\supp (\gp ) := \{
i_1, \ldots , i_m\}$ is called the {\em support} of the idempotent
prime ideal $\gp$. The group $\GL_\infty (\mS_{n-1}, \gp )$ can be
identified with the subgroup $(1+\gp \gp_n)^*$ of the group
$(1+\ga_n)^*$. The group $(1+\gp \gp_n)^*$ contains the descending
chain of normal subgroups
$$(1+\gp\gp_n)^*=(1+\gp\gp_n)^*_1\supset \cdots \supset (1+\gp\gp_n)^*_s\supset
\cdots \supset (1+\gp\gp_n)^*_n=(1+F_n)^*\supset
(1+\gp\gp_n)^*_{n+1}=\{ 1\}$$ where
$(1+\gp\gp_n)^*_s:=(1+\gp\gp_n)^*\cap(1+\ga_{n,s})^*$. Moreover,
the groups $(1+\gp\gp_n)^*_s$ are normal subgroups of the group
$(1+\ga_n)^*$.  The following lemma describes the normal
subgroups $(1+\gp\gp_n)^*_s$.

\begin{lemma}\label{a12May10}
Let $\gp = \gp_{i_1}+\cdots +\gp_{i_m}$ where $i_1, \ldots , i_m$
are distinct elements of the set $\{1, \ldots , n\}$. Then
$$(1+\gp\gp_n)_s^* =\begin{cases}
(1+\sum_{|I|=s, I\in \CJ (\gp )}\gp_I)^*  & \text{if }s=2, \ldots , n-1,\\
(1+F_n)^*& \text{if }s=n, \\
\end{cases}$$
where $\CJ (\gp ) :=\{ J\subseteq \{ 1, \ldots , n\} \, | \, n\in
J, J\cap \supp (\gp ) \neq \emptyset \}$. In particular
$(1+\gp\gp_n)_1^*=(1+\gp\gp_n)_2^*=(1+\gp\gp_n)^*$.
\end{lemma}

{\it Proof}. The case $s=n$ is obvious. So,   we assume that
$s\neq n$. Since the ideals $\gp\gp_n$ and $\ga_{n,s}$ of the
algebra $\mS_n$ are idempotent ideals, $$ \gp \gp_n\cap \ga_{n,s}
= \gp \gp_n\ga_{n,s} = \sum_{\nu =1}^m\gp_{i_\nu}\gp_n\ga_{n,s}=
\sum_{|I|=s, I\in \CJ (\gp )}\gp_I.
$$
Therefore,  $(1+\gp\gp_n)_s^* =(1+\gp\gp_n)^* \cap (1+\ga_{n,s})^*
= (1+\gp\gp_n\cap \ga_{n,s})^* = (1+\sum_{|I|=s, I\in \CJ (\gp
)}\gp_I)^*$.
 $\Box $

$\noindent $

By (\ref{apns1}), there is  a group monomorphism
$$ \v_{n,s}: \frac{(1+\gp\gp_n)^*_s}{(1+\gp\gp_n)^*_{s+1}}\ra
\frac{(1+\ga_{n,s})^*}{(1+\ga_{n,s+1})^*}\ra
\prod_{|I|=s}(1+\bgp_I)^*= \prod_{|I|=s,  I\in \CJ (\gp )
}(1+\bgp_I)^*\times \prod_{|I'|=s, I'\not\in \CJ (\gp
)}(1+\bgp_{I'})^*$$ which is the composition of two group
monomorphisms. By Lemma \ref{a12May10}, 
\begin{equation}\label{Pipns1}
\im (\v_{n,s})\subseteq \prod_{|I|=s, I\in \CJ (\gp
)}(1+\bgp_I)^*.
\end{equation}
For each number $s=2, \ldots , n-1$, consider the following subset
of the group $(1+\gp\gp_n)^*$,
$$\tmE_{n,s}(\gp ):= \prod_{|I|=s, I\in \CJ (\gp )} U_I\ltimes
E_\infty (\mS_{CI}).$$
 It is a product of subgroups of
$(1+\gp\gp_n)_s^*$ in arbitrary order,  but which is assumed to be fixed
for each $s$.

Recall that $(1+\bgp_I)^* = (\mX_{CI}\times U_I)\ltimes E_\infty
(L_{CI})$. Since $\v_{n,s}(\tmE_{n,s}(\gp )(1+\gp
\gp_n)^*_{s+1})=\prod_{|I|=s, I\in \CJ (\gp )} U_I\ltimes E_\infty
(L_{CI})$, we see that there is the group monomorphism
$$\overline{\v}_{n,s}: \frac{(1+\gp\gp_n)^*_s}{ \tmE_{n,s}(\gp )(1+\gp\gp_n)^*_{s+1}}\ra
  \frac{(1+\ga_{n,s})^* }{ \G_{n,s}}\simeq
\mX_{n,s}\simeq \mX_{n,s}'= \prod_{|I|=s, I\in \CJ (\gp
)}\mX_{CI}'\times \prod_{|I'|=s, I'\in \CJ (\gp )}\mX_{CI'}'.$$
Notice that the group $\tmE_{n,s}(\gp )(1+\gp\gp_n)^*_{s+1}$ is a
normal subgroup of $(1+\gp \gp_n)^*_s$. For each number $s=2,
\ldots , n-1$, in the set $\Theta_{n,s}'$ consider  the exact
product of cyclic groups (the order is arbitrary)
\begin{equation}\label{Pttns}
\tTh_{n,s}(\gp )= \tTh_{n,s}^{[1]}(\gp )\times_{ex}
\tTh_{n,s}^{[2]}(\gp ),
\end{equation}
\begin{eqnarray*}
\tTh_{n,s}^{[1]}(\gp )&:= &{}^{exact}\prod_{i\in \supp (\gp )} \;
\prod_{|J'|=s+1, n\in J', J'\cap \supp (\gp )=\{ i\} }\;
\prod_{j'\in J\backslash \{ n, i, m'(J')\}} \langle
\theta_{m'(J'), j'}(J')\rangle ,\\  \tTh_{n,s}^{[2]}(\gp )&:=&
{}^{exact} \prod_{|J|=s+1, n\in J, J\cap \supp (\gp ) \geq 2}\;
\prod_{j\in J\backslash \{ n, m(J)\}}\langle \theta_{m(J),
j}(J)\rangle ,
\end{eqnarray*}
where $m'(J')$ is the largest element of the set $J'\backslash \{
n, i\}$ and $m(J)$ is the largest  element of the set $J\backslash
n$. Notice that $\tTh_{n,2}(\gp )=\tTh_{n,2}^{[2]}(\gp )$ as the
set $\tTh_{n,2}^{[1]}(\gp )$ is an empty set.

 By (\ref{Pipns1}),  $\im
(\overline{\v}_{n,s})\subseteq \prod_{|I|=s,  I\in \CJ (\gp
)}\mX_{CI}'$ and
\begin{eqnarray*}
\prod_{|I|=s,  I\in \CJ (\gp )}\mX_{CI}' & =&\prod_{|I|=s,  I\in
\CJ (\gp )}\; \prod_{i\in CI}\Z (i,I)= \prod_{|J|=s+1, n\in J,
J\cap \supp (\gp )\neq \emptyset}\;\;
\prod_{j\in J\backslash n, (J\backslash n )\cap \supp (\gp )\neq \emptyset}\Z (j, J\backslash j) \\
 &=& \prod_{i\in \supp (\gp )}
\;\prod_{|J'|=s+1, n\in J', J'\cap \supp (\gp )=\{ i\} }\;
\prod_{j'\in J'\backslash \{ n, i\}} \Z(j',J'\backslash
j')\\
&\times & \prod_{|J|=s+1, n\in J, |J\cap \supp (\gp )| \geq 2}\;
\prod_{j\in J\backslash  n} \Z (j,J\backslash j).
\end{eqnarray*}
 Recall that $\im (\psi_{n,s}')= \psi_{n,s}'(\Theta_{n,s}) = \mK_{n,s}=\bigcap_{|J|=s+1}\ker (\chi_J')$,
  by Theorem \ref{24May9} and (\ref{ptijJ1}). The following
  argument is the key moment in the proof of Theorem
  \ref{12May10}. For each number $s=2, \ldots , n-1$,
\begin{eqnarray*}
 \im (\overline{\v}_{n,s})&\subseteq & \im (\psi_{n,s}')\bigcap \prod_{|I|=s,
 I\in \CJ (\gp )}\mX_{CI}'=\bigcap_{|I|=s+1}\ker (\chi_J')\bigcap \prod_{|I|=s,
 I\in \CJ (\gp )}\mX_{CI}'\\
 &=&\prod_{i\in \supp (\gp )}
\;\prod_{|J'|=s+1, n\in J', J'\cap \supp (\gp )=\{ i\} }\;
\prod_{j'\in J\backslash \{ n, i, m'(J')\}} \Z (-(m'(J'),
J'\backslash
m'(J'))+(j', J\backslash j'))\\
&\times & \prod_{|J|=s+1, n\in J, |J\cap \supp (\gp )| \geq 2}\;
\prod_{j\in J\backslash \{ n, m(J)\}} \Z (-(m(J), J\backslash
m(J))+(j, J\backslash j))\\
&\stackrel{{\rm by}\; (\ref{ptijJ})}{=}&
\psi_{n,s}'(\tTh_{n,s}(\gp ))= \overline{\v}_{n,s}'(\tTh_{n,s}(\gp
)\tmE_{n,s}(\gp )(1+\gp \gp_n)^*_{s+1}).
\end{eqnarray*}
The first equality above follows from the decomposition of  the
abelian group $\prod_{|I|=s,  I\in \CJ (\gp )}\mX_{CI}'$ above and
the definition of the homomorphisms $\chi_J'$.  It follows that
\begin{equation}\label{Ppns1}
(1+\gp\gp_n)^*_s= \tTh_{n,s}(\gp )\times_{ex}\tmE_{n,s}(\gp
)(1+\gp\gp_n)^*_{s+1}, \; \; s=2, \ldots , n-1.
\end{equation}

\begin{theorem}\label{12May10}
Let $\gp$ be a nonzero idempotent prime ideal of the algebra
$\mS_{n-1}$. Then the group $\GL_\infty (\mS_{n-1}, \gp ) =
(1+\gp\gp_n)^* $ is equal to $ \tTh_{n,2}(\gp )\tmE_{n,2}(\gp
)\cdots \tTh_{n,n-1}(\gp )\tmE_{n,n-1}(\gp )$.
 Moreover,
$$ (1+\gp\gp_n)_s^*=\begin{cases}
\tTh_{n,2}(\gp )\tmE_{n,2}(\gp )\cdots \tTh_{n,n-1}(\gp )\tmE_{n,n-1}(\gp )& \text{if }s=1,\\
\tTh_{n,s}(\gp )\tmE_{n,s}(\gp )\cdots \tTh_{n,n-1}(\gp )\tmE_{n,n-1}(\gp )& \text{if }s=2, \ldots, n-1,\\
(1+F_n)^* &  \text{if }s=n.\\
\end{cases}$$
\end{theorem}

{\it Proof}. By Proposition 3.10, \cite{Snaut}, we have the
inclusion  $(1+\gp \gp_n)_n^*=(1+F_n)^* \subseteq \tmE_{n,n-1}(\gp
)$. Now, the theorem follows from (\ref{Ppns1}). $\Box $


\begin{lemma}\label{a13May10}
Let $\mS_1(\L )$ be the algebra $\mS_1$ over the ring $\L$ from
Lemma \ref{a8May10}. Then,  for each element $\l \in D\backslash
\{ -1\}$,
\begin{equation}\label{CcomE}
\Bigl(
\begin{matrix} 1  & 0\\ -\frac{ey}{1+\l e} & 1
\end{matrix}\Bigr)\,
\Bigl(
\begin{matrix} 1 & \l ex \\  0& 1
\end{matrix}\Bigr)\,
\Bigl(
\begin{matrix} 1 &0 \\ey  & 1
\end{matrix}\Bigr)\,
\Bigl(
\begin{matrix} 1 & -\frac{\l e}{1+\l e} x \\  0& 1
\end{matrix}\Bigr)\,
= \Bigl(
\begin{matrix} 1+\l e& 0\\ 0 & \frac{1}{1+\l e}
\end{matrix}\Bigr)\,
\Bigl(
\begin{matrix} 1-\frac{\l }{1+\l}eE_{00} & 0\\ 0 &1
\end{matrix}\Bigr).
\end{equation}
where $E_{00}:=1-xy$ and  $\frac{1}{1+\l e}=1-\frac{\l}{1+\l} e$,
by Lemma \ref{a8May10}.(2).
\end{lemma}

{\it Proof}. The RHS of the equality (\ref{CcomE}) is the product
of four matrices, say $A_1\cdots A_4$.
$$ A_1A_2A_3=\Bigl(
\begin{matrix} 1  & \l ex\\ -\frac{ey}{1+\l e} & \frac{1}{1+\l e}
\end{matrix}\Bigr)\,
\Bigl(
\begin{matrix} 1 &0 \\ey  & 1
\end{matrix}\Bigr)=\Bigl(
\begin{matrix} 1+\l exy  & \l ex\\ 0 & \frac{1}{1+\l e}
\end{matrix}\Bigr), \, \;\;\; A_1\cdots A_4=\Bigl(
\begin{matrix} 1+\l exy  & 0\\ 0 & \frac{1}{1+\l e}
\end{matrix}\Bigr), $$
since $(1+\l exy) (-\frac{\l e}{1+\l e} x) +\l ex = -\frac{\l
e}{1+\l e}(1+\l e) x+\l ex= 0$. The product $A_1\cdots A_4$
coincides with the product `$A_1\cdots A_4$' in the proof of Lemma
\ref{b8May10}, and so the equality (\ref{CcomE}) follows  from
(\ref{ABcomE}).  $\Box $


\begin{lemma}\label{b13May10}
$E_\infty' (\mS_{n-1}, \gp ) \supseteq \tmE_{n,s}$ for all $s=3,
\ldots , n-1$ and $E_\infty' (\mS_{n-1}, \gp ) \supseteq E_\infty
(\mS_{CI})$ for all sets $I\in \CJ (\gp )$ such that $|I|=2$.
\end{lemma}

{\it Proof}. We have to show that the group $E_\infty' (\mS_{n-1},
\gp )$ contains the groups  $E_\infty (\mS_{CI})$ for all subsets
$I\in \CJ (\gp )$ such that $|I|=2, \ldots , n-1$, and the groups
$U_I$ for all subsets  $I\in \CJ (\gp )$ such that $|I|=3, \ldots
, n-1$. By Lemma \ref{a13May10}, the groups $U_I$ belong to the
group $E_\infty' (\mS_{n-1}, \gp )$. Indeed, by (\ref{muIG}), each
element of the group $U_I$ is a matrix $u = \Bigl(
\begin{matrix} 1+\mu e_{I\backslash n}  &0 \\0 &1
\end{matrix}\Bigr)$ for some scalar $\mu \in K\backslash \{ -1\}$.
 Since $I\in \CJ (\gp )$ and $|I|\geq 3$, we can choose a number
 $j\in I\backslash n$ such that $(I\backslash \{ j,n\} ) \cap
 \supp (\gp )\neq 0$. Then $e_{I\backslash n}= e\cdot E_{00}(j)$
 where $e=e_{I\backslash \{ j,n\}}\in \gp$. By Lemma
 \ref{a13May10}, the matrix $u$ belongs to the group $E_\infty
 (\mS_{n-1}, \gp )$, since the map $\v : K\backslash \{ -1\} \ra K\backslash \{
 -1\}$, $\l \mapsto -\frac{\l}{1+\l}$, is a bijection.

 The group $E_\infty (\mS_{CI})$ is generated by the
elementary matrices $u=1+aE_{\alpha\beta}(I)$ where $a\in
\mS_{CI}$, $\alpha = (\alpha_i)_{i\in I}$, $\beta =
(\beta_i)_{i\in I}\in \N^I$ and $\alpha \neq \beta$. If
$\alpha_n\neq \beta_n$ then $u=1+(a\prod_{i\in I, i\neq
n}E_{\alpha_i\beta_i}(i))E_{\alpha_n\beta_n}(n)\in E_\infty'
(\mS_{n-1}, \gp )$, since $I\in \CJ (\gp )$. If $\alpha_n=\beta_n$
then choose an element $\g \in \N^I$ such that $\g_n\neq
\alpha_n$, and so $\g \neq \alpha$ and $\g \neq \beta$. Since the
elements $1+E_{\alpha\g}$ and $1+aE_{\g\beta}$ belong to the group
$E_\infty' (\mS_{n-1}, \gp )$ (by the previous case), so does their
group commutator
$$E_\infty (\mS_{n-1}, \gp )\ni [1+E_{\alpha\g}, 1+aE_{\g\beta}]=1+aE_{\alpha\beta}=u.$$
Therefore, $E_\infty (\mS_{CI})\subseteq E_\infty' (\mS_{n-1}, \gp
)$. $\Box$

\begin{lemma}\label{a14May10}
Let $J=\{ i,j,n\}$ where the numbers $i$, $j$ and $n$ are
distinct. Let $I=\{ k,n\}$ where $k\neq n$, and $\l \in K^*$. Then
$$ [\th_{ij}(J), \mu_I(\l )]= \begin{cases}
1 & \text{if }k\neq i, k\neq j,\\
1+(\l^{-1}-1)e_J= \mu_J(\l )^{-1}& \text{if }k=i, \\
1+(\l -1)E_{11}(j)e_ie_n& \text{if }k=j.\\
\end{cases}
$$
\end{lemma}

{\it Proof}. Let $c$ be the group commutator, $J'=\{ i,j\}$,
$\th_{ij}=\th_{ij}(J)$ and $\th_{ij}'=\th_{ij}(J')$. Since
$\th_{ij}^{\pm 1}e_n = e_n\th_{ij}^{\pm 1}= \th_{ij}'^{\pm 1}e_n =
e_n\th_{ij}'^{\pm 1}$ and $\th_{ij}'^{-1}=\th_{ji}'$, we see that
$$ c= \th_{ij}(1+(\l -1)e_ke_n) \th_{ij}^{-1} \mu_I^{-1}(\l ) = (1+(\l
-1)\th_{ij}'e_k\th_{ji}'e_n)\mu_I^{-1}(\l ).$$ If $k\neq i$ and
$k\neq j$ then  the elements $\th_{ij}'$ and $e_k$ commute and we
get $c=\mu_I(\l ) \mu_I(\l )^{-1}=1$.

If $k=i$ then $\th_{ij}'e_i=x_je_i$ and $e_i \th_{ji}'=e_iy_j$, by
(\ref{t12s}), and so
\begin{eqnarray*}
  c&=& (1+(\l -1) x_jy_je_ie_n)\mu_I(\l )^{-1} = (\mu_I(\l ) - (\l
-1)e_J) \mu_I(\l )^{-1} \\
 &=&1-(\l -1)e_J(1+(\l^{-1}-1)e_I)= 1-\frac{\l -1}{\l}e_J=
 1+(\l^{-1}-1)e_J = \mu_J(\l )^{-1}.
\end{eqnarray*}
If $k=j$ then $\th_{ij}'e_j= y_ie_j+E_{10}(j)e_i$ and $ e_j
\th_{ji}'=x_ie_j+E_{01}(j)e_i$, by (\ref{t12s}), and so
\begin{eqnarray*}
  c&=& (1+(\l -1)(y_ie_j+E_{10}(j)e_i) (x_ie_j+E_{01}(j)e_i)e_n)
  \mu_I(\l )^{-1}\\
  & =& (\mu_I(\l ) + (\l
-1)E_{11}(j)e_ie_n) \mu_I(\l )^{-1}
 =1+(\l -1)E_{11}(j)e_ie_n. \;\; \Box
\end{eqnarray*}


Let $A$ and $B$ be subgroups/subsets of a group $G$. The  {\em
commutant} $[A,B]$ is the subgroup of $G$ generated by all the
group commutators $[a,b]=aba^{-1}b^{-1}$ where $a\in A$ and  $b\in
B$. For an element $g\in G$, let  $\o_g : x\mapsto gxg^{-1}$ be
the inner automorphism of the group $G$ determined by the element
$g$. We can easily verify that for all elements $a_1, a_2, b_1,
b_2\in G$, 
\begin{equation}\label{a1b1}
 [a_1b_1, a_2b_2]= \o_{a_1}([b_1, a_2]) \o_{a_1a_2}([b_1, b_2]) [
a_1, a_2] \o_{a_2}([a_1, b_2]).
\end{equation}

{\bf The normal subgroup $\CE (\mS_{n-1}, \gp )$}. Consider the
subgroup
$$\CE (\mS_{n-1}, \gp ):= \prod_{|I|=2, I\in \CJ (\gp )}
E_\infty (\mS_{CI})\cdot (1+\gp\gp_n)^*_3$$ of the group
$(1+\gp\gp_n)^*=\GL_\infty (\mS_{n-1}, \gp )$. By (\ref{Ppns1}),
the group $(1+\gp\gp_n)^*$ is the exact product of sets,
\begin{equation}\label{ppex1}
(1+\gp\gp_n)^*=\tTh_{n,2}(\gp ) \times_{ex}{}^{exact}\prod_{|I|=2,
I\in \CJ (\gp )}U_I\times_{ex}\CE (\mS_{n-1},\gp ).
\end{equation}
By the very definition, the subgroup $\CE (\mS_{n-1}, \gp )$ is a
{\em normal} subgroup of $(1+\gp\gp_n)^*$ (see the definition of
the map $\v_{n,s}$). There is the inclusion 
\begin{equation}\label{ppex2}
[\tTh_{n,2}(\gp ) , \tTh_{n,2}(\gp ) ] \subseteq (1+\gp\gp_n)^*_3
\end{equation}
which is obvious due to the fact that the image of each element
$\th_{ij}(J)$ (where $|J|=3$  and $J\in \CJ (\gp )$) under the map
$\v_{n,s}$ is the direct product of two `diagonal' matrices with
entries in (commutative) Laurent polynomial algebras, hence all
the images commute.
\begin{theorem}\label{A13May10}
$\CE (\mS_{n-1}, \gp ) = E_\infty (\mS_{n-1}, \gp )=E_\infty' (\mS_{n-1}, \gp )$.
\end{theorem}
{\it Proof}. Recall that $\GL_\infty (R, \ga) /E_\infty (R, \ga )$ is an abelian group  for any ring $R$ and  ideal
$\ga$ of $R$, \cite{Bass-book-K-theory}. By (\ref{ppex1}), Lemma
\ref{a14May10} and (\ref{ppex2}), the factor group
$(1+\gp\gp_n)^*/\CE (\mS_{n-1}, \gp )$ is abelian.

 Let us show that $E_\infty'(\mS_{n-1}, \gp )  \subseteq \CE := \CE (\mS_{n-1}, \gp )$. We have to show that $1+\gp E_{ij}(n)\subseteq  \CS$ for all $i\neq j$. Since
 $$ 1+\gp E_{ij}(n)=1+(\gp_{i_1}+\cdots +\gp_{i_m}) E_{ij}(n)=\prod_{\nu =1}^m(1+\gp_{i_\nu} E_{ij}(n)),$$it suffices to show that $1+\gp_{i_\nu} E_{ij}(n)\subseteq \CE$ for all $\nu = 1, \ldots , m$ and $i\neq j$, but this is obvious since
 $$1+\gp_{i_\nu} E_{ij}(n)\subseteq E_\infty (\mS_{CI})\subseteq \CE$$ where $I=\{ i_\nu, n \}$ (and so $|I|=2$ and $I\in \CJ (\gp )$), see the definition of $\CE$.


To finish the proof of the theorem it suffices to show that
 $\CE
(\mS_{n-1}, \gp ) \subseteq  E_\infty' (\mS_{n-1}, \gp )$ (since then the group $E_\infty' (\mS_{n-1}, \gp )$ is a normal subgroup of $\GL_\infty (\mS_{n-1}, \gp )$ as $\GL_\infty (\mS_{n-1}, \gp ) / \CE$ is an abelian group  and $\CE\subseteq E_\infty'(\mS_{n-1},\gp )$; hence $E_\infty (\mS_{n-1}, \gp ) = E_\infty' (\mS_{n-1}, \gp )$ and $E_\infty' (\mS_{n-1}, \gp )=\CE$). By
Theorem \ref{12May10},
$$\CE
(\mS_{n-1}, \gp ) = \prod_{|I|=2, I\in \CJ (\gp )}E_\infty
(\mS_{CI})\cdot \tTh_{n,3}(\gp ) \tmE_{n, 3}(\gp ) \cdots
\tTh_{n,n-1}(\gp ) \tmE_{n, n-1}(\gp ).$$ By Lemma \ref{b13May10},
the inclusion $\CE (\mS_{n-1}, \gp ) \subseteq E_\infty'
(\mS_{n-1}, \gp )$ holds iff $\tTh_{n,s}(\gp ) \subseteq E_\infty'
(\mS_{n-1}, \gp )$ for all $s=3, \ldots , n-1$ iff
$\tTh_{n,s}^{[1]}(\gp ), \tTh_{n,s}^{[2]}(\gp ) \subseteq E_\infty'
(\mS_{n-1}, \gp )$ for all $s=3, \ldots , n-1$.

Fix an element $\th$ such that either $\th \in
\tTh_{n,s}^{[1]}(\gp )$ or $\th \in \tTh_{n,s}^{[2]}(\gp )$, i.e.
either $\th = \th_{m'(J'), j'}(J')$ or $\th = \th_{m (J), j}(J)$,
see (\ref{Pttns}). In the second case, without less of generality
we may assume that $m(J) \not\in J\cap \supp (\gp )$, by changing,
if necessary, the order in the set $J$ (or simply by taking a
suitable element). In both cases, we can choose an element, say
$k\in J\cap \supp (\gp )$, such that $k\not\in \{ m'(J'), j'\}$,
in the first case, and $k\not\in \{ m(J), j\}$, in the second
case. In both cases, we can write $\th = \th_{ij}(J)$ where $k\in
J\cap \supp (\gp )$ and $k\not\in \{ i,j\}$. As we have seen in
Section \ref{EEE},
$$\th_{ij}(J) = e_k\th_{ij}(J\backslash k) +1-e_k.$$ By Theorem
\ref{17Apr10}, $\th_{ij}(J\backslash k) \in E_\infty
(\bigotimes_{l=1, l\neq k}^{n-1}\mS_1(l))\subseteq \GL_\infty
(\bigotimes_{l=1, l\neq k}^{n-1}\mS_1(l))$. Under the algebra
monomorphism
$$\GL_\infty (\bigotimes_{l=1, l\neq
k}^{n-1}\mS_1(l))\ra \GL_\infty (\mS_{n-1}, \gp ), \;\; a\mapsto
e_ka+1-e_k,$$ the group of elementary matrices $E_\infty
(\bigotimes_{l=1, l\neq k}^{n-1}\mS_1(l))$ is mapped into the
group of $\gp$-elementary matrices $E_\infty' (\mS_{n-1}, \gp )$
since $e_k\in \gp$. Therefore, $\th \in E_\infty' (\mS_{n-1}, \gp
)$. The proof of the theorem is complete. $\Box $

\begin{theorem}\label{13May10}
Let $\gp$ be a nonzero idempotent prime ideal of the algebra
$\mS_{n-1}$ and $m={\rm ht} (\gp )$ be its height. Then (below is the
direct product of groups)
$$ {\rm K}_1(\mS_{n-1}, \gp ) \simeq \prod_{ \{i>j \, | \, i,j\in
\supp (\gp )\}} \langle \th_{ij}(\{ i,j,n\})\rangle \times
\prod_{k\in \supp (\gp )} U_{ \{ k,n\} } \simeq \begin{cases}
K^*, & \text{if }m=1,\\
\Z^{m\choose 2}\times K^{*m}& \text{if }m> 1.\\
\end{cases}$$
The group $\GL_\infty (\mS_{n-1}, \gp )$ is generated by the
elements $\th_{ij}:=\th_{ij}(\{ i,j,n\})$  (where $i>j$ and
$i,j\in \supp (\gp )$) and the groups $E_\infty (\mS_{n-1}, \gp
)$, $U_{\{ k,n\} }$ where $k\in \supp (\gp )$. Moreover, each
element $a$ of the group $\GL_\infty (\mS_{n-1}, \gp )$ is the
unique product (the order is arbitrary) 
\begin{equation}\label{aupSp}
a=\prod_{ \{i>j \, | \, i,j\in \supp (\gp )\}}
\th_{ij}^{n_{ij}}\cdot \prod_{k\in \supp (\gp )} \mu_{ \{ k,n\}
}(\l_k) \cdot e
\end{equation}
where $n_{ij}\in \Z$, $\l_k\in K^*$ and $e\in E_\infty (\mS_{n-1},
\gp )$.
\end{theorem}

{\it Proof}. The theorem follows from the equality (\ref{ppex1})
and Theorem \ref{A13May10}. $\Box $

$\noindent $

We can find effectively (in finitely many steps) the decomposition
(\ref{aupSp}) (Corollary \ref{a15May10}). For, we introduce
several explicit group homomorphisms.

$\noindent $

{\it Definition}. For each nonempty subset $I$ of $\{ 1, \ldots ,
n\}$ with $s=|I|<n$ and for each element $j\in CI$, define the
group homomorphism $ \det_I:(1+\ga_{n,s})^*\ra L^*_{CI}$ as the
composition of the group homomorphisms (see (\ref{apns1}))
$$(1+\ga_{n,s})^*\stackrel{\psi_{n,s}}{\ra}\prod_{|J|=s}(1+\bgp_J)^*\stackrel{{\rm pr}_I}{\ra}(1+\bgp_I)^*\simeq
\GL_\infty (L_{CI})\stackrel{\det}{\ra} L_{CI}^*$$ where ${\rm
pr}_I$ is the projection map. Define the group homomorphism $
\deg_{n,I,j}:(1+\ga_{n,s})^*\ra \Z$ as the composition of the
group homomorphisms  $(1+\ga_{n,s})^*\stackrel{\det_I}{\ra}
L_{CI}^*\stackrel{\deg_{x_j}}{\ra}\Z$ where  $\deg_{x_j}$ is the
degree in $x_j$ of monomial ($\deg_{x_j}(\l \prod_{i\in
CI}x_i^{\alpha_i}) = \alpha_j$ where $\l \in K^*$ and $\alpha_i\in
\Z$).

\begin{lemma}\label{a18Apr10}
Let $n\geq 3$ and $s=1, \ldots , n-1$. Then for all subsets $I$
and $J$ of the set $\{ 1, \ldots , n\}$  such that $ |I|=s$, $|J|=
s+1$ and $n\in J$,
$$ \deg_{n,I, i} (\th_{m(J), j}(J))=
\begin{cases}
-1& \text{if }I=J\backslash m(J), i=m(J),\\
1& \text{if }I=J\backslash j, i=j,\\
0& \text{otherwise}.\\
\end{cases}$$
where $i\in CI$ and $j\in J\backslash m(J)$.

\end{lemma}

{\it Proof}. The result follows at once from the equality
$\th_{m(J), j}=(1+(y_{m(J)}-1)e_{J\backslash
m(J)})(1+(x_j-1)e_{J\backslash j})$. $\Box $

\begin{corollary}\label{a15May10}
Given  a product decomposition   (\ref{aupSp}) for an element $a\in
\GL_\infty (\mS_{n-1}, \gp )$, we have
\begin{eqnarray*}
n_{ij}& = & \deg_{n, \{ i,n\}, j}(a),\\
\l_k &=& \det_{\{k,n\}}(a\cdot \prod_{ \{i>j \, | \, i,j\in \supp
(\gp )\}} \th_{ij}^{-n_{ij}}),\\
e&=&(\prod_{ \{i>j \, | \, i,j\in \supp (\gp )\}}
\th_{ij}^{n_{ij}}\cdot  \prod_{k\in \supp (\gp )} \mu_{ \{ k,n\}
}(\l_k))^{-1} a.
\end{eqnarray*}
\end{corollary}
{\it Proof}. By Lemma \ref{a18Apr10}, $\deg_{n, \{ i,n\}, j}(a)=
n_{ij}\deg_{n, \{ i,n\}, j}(\th_{ij})=n_{ij}$. Similarly,
$$\det_{\{k,n\}}(a\cdot \prod_{ \{i>j \, | \, i,j\in \supp (\gp
)\}} \th_{ij}^{-n_{ij}}) = \det_{\{k,n\}}(\mu_{\{ k,n\} }
(\l_k))=\l_k.$$ The rest is obvious. $\Box $

$\noindent $

Corollary \ref{a15May10} gives an effective criterion of whether
an element $a\in \GL_\infty (\mS_{n-1}, \gp )$ is a product of
$\gp$-elementary matrices.
\begin{corollary}\label{b15May10}
Let $a\in \GL_\infty (\mS_{n-1}, \gp )$. Then $a\in E_\infty
(\mS_{n-1}, \gp )$ iff all $n_{ij}=$ and $\l_k=1$ iff $\deg_{n, \{
i,n\}, j}(a)=1$ for all $i>j$ such that $i,j\in \supp (\gp )$, and
$\det_{\{k,n\}}(a)=1$ for all $k\in \supp (\gp )$.
\end{corollary}

$${\bf Acknowledgements}$$

$\noindent $

I would like to thank A. Bak and  C. A. Weibel  for the comments.

Department of Pure Mathematics

University of Sheffield

Hicks Building

Sheffield S3 7RH

UK

email: v.bavula@sheffield.ac.uk

\end{document}